\documentclass[12pt]{article}
\usepackage{mathrsfs}
\usepackage{amssymb}
\usepackage{amsfonts}
\usepackage{latexsym}
\usepackage{amsthm}
\usepackage{pictex}

\def\!{\mskip-\thinmuskip}

\usepackage{graphicx}

\newcommand{\er}[1]{{\rm(\ref{#1})}}
\def\lb{\label}

\theoremstyle{plain}
\newtheorem{theorem}{\bf Theorem}[section]
\newtheorem{lemma}[theorem]{\bf Lemma}
\newtheorem{corollary}[theorem]{\bf Corollary}

\theoremstyle{remark}

\begin{document}

\def\a{\alpha} 
\def\b{\beta}  
\def\g{\gamma} \def\G{\Gamma}
\def\d{\delta} \def\D{\Delta}
\def\c{\chi}
\def\z{\zeta}
\def\e{\eta}
\def\f{\phi}     \def\F{\Phi} 
\def\k{\kappa} 
\def\l{\lambda} \def\L{\Lambda}
\def\m{\mu}
\def\n{\nu}  
\def\o{\omega} \def\O{\Omega}
\def\p{\psi}   \def\P{\Psi} 
\def\r{\rho}  
\def\s{\sigma}  \def\S{\Sigma} 
\def\vT{\Theta}
\def\s{\sigma} \def\S{\Sigma} 
\def\x{\xi}   \def\X{\Xi}
\def\ve{\varepsilon}
\def\vt{\vartheta}
\def\vp{\varphi}
\def\vk{\varkappa}
\def\t{\tau}

\def\cA{{\cal A}} \def\bA{{\bf A}}  \def\mA{{\mathscr A}}
\def\cB{{\cal B}} \def\bB{{\bf B}}  \def\mB{{\mathscr B}}
\def\cC{{\cal C}} \def\bC{{\bf C}}  \def\mC{{\mathscr C}}
\def\cD{{\cal D}} \def\bD{{\bf D}}  \def\mD{{\mathscr D}}
\def\cE{{\cal E}} \def\bE{{\bf E}}  \def\mE{{\mathscr E}}
\def\cF{{\cal F}} \def\bF{{\bf F}}  \def\mF{{\mathscr F}}
\def\cG{{\cal G}} \def\bG{{\bf G}}  \def\mG{{\mathscr G}}
\def\cH{{\cal H}} \def\bH{{\bf H}}  \def\mH{{\mathscr H}}
\def\cI{{\cal I}} \def\bI{{\bf I}}  \def\mI{{\mathscr I}}
\def\cJ{{\cal J}} \def\bJ{{\bf J}}  \def\mJ{{\mathscr J}}
\def\cK{{\cal K}} \def\bK{{\bf K}}  \def\mK{{\mathscr K}}
\def\cL{{\cal L}} \def\bL{{\bf L}}  \def\mL{{\mathscr L}}
\def\cM{{\cal M}} \def\bM{{\bf M}}  \def\mM{{\mathscr M}}
\def\cN{{\cal N}} \def\bN{{\bf N}}  \def\mN{{\mathscr N}}
\def\cO{{\cal O}} \def\bO{{\bf O}}  \def\mO{{\mathscr O}}
\def\cP{{\cal P}} \def\bP{{\bf P}}  \def\mP{{\mathscr P}}
\def\cQ{{\cal Q}} \def\bQ{{\bf Q}}  \def\mQ{{\mathscr Q}}
\def\cR{{\cal R}} \def\bR{{\bf R}}  \def\mR{{\mathscr R}}
\def\cS{{\cal S}} \def\bS{{\bf S}}  \def\mS{{\mathscr S}}
\def\cT{{\cal T}} \def\bT{{\bf T}}  \def\mT{{\mathscr T}}
\def\cU{{\cal U}} \def\bU{{\bf U}}  \def\mU{{\mathscr U}}
\def\cV{{\cal V}} \def\bV{{\bf V}}  \def\mV{{\mathscr V}}
\def\cW{{\cal W}} \def\bW{{\bf W}}  \def\mW{{\mathscr W}}
\def\cX{{\cal X}} \def\bX{{\bf X}}  \def\mX{{\mathscr X}}
\def\cY{{\cal Y}} \def\bY{{\bf Y}}  \def\mY{{\mathscr Y}}
\def\cZ{{\cal Z}} \def\bZ{{\bf Z}}  \def\mZ{{\mathscr Z}}

\def\Z{{\Bbb Z}}
\def\R{{\Bbb R}}
\def\C{{\Bbb C}}
\def\T{{\Bbb T}}
\def\N{{\Bbb N}}
\def\S{{\Bbb S}}
\def\H{{\Bbb H}}
\def\dD{{\Bbb D}}

\def\ma{\left(\begin{array}{cc}}    \def\am{\end{array}\right)}
\def\iint{\int\!\!\!\int}
\def\lt{\biggl}                     \def\rt{\biggr}
\let\ge\geqslant                   \let\le\leqslant
\def\[{\begin{equation}}            \def\]{\end{equation}}
\def\wt{\widetilde}                 \def\pa{\partial}
\def\sm{\setminus}                  \def\es{\emptyset}
\def\no{\noindent}                  \def\ol{\overline}
\def\iy{\infty}                     \def\ev{\equiv}
\def\/{\over}
\def\we{\wedge}
\def\ts{\times}
\def\os{\oplus}
\def\ss{\subset}
\def\h{\hat}
\def\wh{\widehat}
\def\Ra{\Rightarrow}
\def\ra{\rightarrow}
\def\la{\leftarrow}
\def\da{\downarrow}
\def\ua{\uparrow}
\def\lra{\leftrightarrow}
\def\Lra{\Leftrightarrow}
\def\Re{\mathop{\rm Re}\nolimits}
\def\Im{\mathop{\rm Im}\nolimits}
\def\supp{\mathop{\rm supp}\nolimits}
\def\sign{\mathop{\rm sign}\nolimits}
\def\Ran{\mathop{\rm Ran}\nolimits}
\def\Ker{\mathop{\rm Ker}\nolimits}
\def\Tr{\mathop{\rm Tr}\nolimits}
\def\const{\mathop{\rm const}\nolimits}
\def\Wr{\mathop{\rm Wr}\nolimits}
\def\diag{\mathop{\rm diag}\nolimits}
\def\dist{\mathop{\rm dist}\nolimits}

\def\th{\theta}
\def\dlint{\displaystyle\int\limits}
\def\iintt{\mathop{\int\!\!\int\!\!\dots\!\!\int}\limits}
\def\intt{\mathop{\int\int}\limits}
\def\lim{\mathop{\rm lim}\limits}
\def\mult{\!\cdot\!}
\def\BBox{\hspace{1mm}\vrule height6pt width5.5pt depth0pt \hspace{6pt}}
\def\1{1\!\!1}
\newcommand{\bwt}[1]{{\mathop{#1}\limits^{{}_{\,\bf{\sim}}}}\vphantom{#1}}
\newcommand{\bhat}[1]{{\mathop{#1}\limits^{{}_{\,\bf{\wedge}}}}\vphantom{#1}}
\newcommand{\bcheck}[1]{{\mathop{#1}\limits^{{}_{\,\bf{\vee}}}}\vphantom{#1}}
\def\nh{\bhat}
\def\nc{\bcheck}
\newcommand{\oo}[1]{{\mathop{#1}\limits^{\,\circ}}\vphantom{#1}}
\newcommand{\po}[1]{{\mathop{#1}\limits^{\phantom{\circ}}}\vphantom{#1}}
\def\ctg{\mathop{\rm ctg}\nolimits}
\def\notto{\to\!\!\!\!\!\!\!/\,\,\,}

\def\pgbrk{\pagebreak}
\def\Dis{\mathop{\rm Dis}\nolimits}

\def\Twelve{
\font\Tenmsa=msam10 scaled 1200
\font\Sevenmsa=msam7 scaled 1200
\font\Fivemsa=msam5 scaled 1200
\textfont\msbfam=\Tenmsb
\scriptfont\msbfam=\Sevenmsb
\scriptscriptfont\msbfam=\Fivemsb

\font\Teneufm=eufm10 scaled 1200
\font\Seveneufm=eufm7 scaled 1200
\font\Fiveeufm=eufm5 scaled 1200
\textfont\eufmfam=\Teneufm
\scriptfont\eufmfam=\Seveneufm
\scriptscriptfont\eufmfam=\Fiveeufm}

\def\Ten{
\textfont\msafam=\tenmsa
\scriptfont\msafam=\sevenmsa
\scriptscriptfont\msafam=\fivemsa

\textfont\msbfam=\tenmsb
\scriptfont\msbfam=\sevenmsb
\scriptscriptfont\msbfam=\fivemsb

\textfont\eufmfam=\teneufm
\scriptfont\eufmfam=\seveneufm
\scriptscriptfont\eufmfam=\fiveeufm}

\title { Spectral estimates for Schr\"odinger
operator with periodic matrix potentials on the real line}

\author{Dmitri Chelkak
\begin{footnote}
{Dept. of Math. Analysis, Math. Mech. Faculty,
St.Petersburg State University. Universitetskij pr. 28, Staryj Petergof, 198504 St.Petersburg,
Russia, e-mail: delta4@math.spbu.ru}
\end{footnote}${}^{\!\!\!}$
and Evgeny Korotyaev
\begin{footnote}
{Correspondence author. Institut f\"ur  Mathematik,  Humboldt Universit\"at zu Berlin,
Rudower Chaussee 25, 12489, Berlin, Germany, e-mail:
evgeny@math.hu-berlin.de}
\end{footnote}
}

\maketitle

\begin{abstract}
\no We consider the Schr\"odinger operator on the real line with a
$N\ts N$ matrix  valued periodic potential, $N>1$. The spectrum
of this operator is absolutely continuous and consists of
intervals separated by gaps. We define the Lyapunov function,
which is analytic on an associated N-sheeted Riemann surface. On each sheet
the Lyapunov function has the standard properties of the Lyapunov
function for the scalar case. The Lyapunov function has (real or complex) branch points, which we call resonances.
We determine the asymptotics of the periodic, anti-periodic spectrum and of the
resonances at high energy (in terms of the Fourier coefficients of the potential). We show that there exist two types of
gaps: i) stable gaps, i.e., the endpoints are periodic and
anti-periodic eigenvalues, ii) unstable (resonance) gaps, i.e., the
endpoints are resonances (real branch points). Moreover, the
following results are obtained: 1) we define the quasimomentum as
an analytic function on the Riemann surface of the Lyapunov function; various properties
and estimates of the quasimomentum are obtained, 2) we construct
the conformal mapping with real part given by the 
 integrated density of states and imaginary part given by the
 Lyapunov exponent. We obtain various properties of this
 conformal mapping, which are similar to the case N=1,
3) we determine various new trace formulae for potentials, the
integrated density of states  and the
 Lyapunov exponent, 4) a priori estimates of gap
lengths in terms of  potentials are obtained.

\end{abstract}


\section {Introduction and main results}
\setcounter{equation}{0}

We consider the self-adjoint operator $\mL y=-y''+V(t)y,$ acting
in $L^2(\R)^N, N\ge 2$ where the symmetric 1-periodic $N\ts N$ matrix
potential $V$ belongs to the real Hilbert space $\mH$ given by
$$
\mH=\lt\{V=V^*=\{V_{jk}(t)\}_{j,k=1}^N,\ \ t\in \R/\Z,\ \ \
\|V\|^2=\int_0^1\Tr V^2(t)dt<\iy\rt\}.
$$
It is well known (see [DS]) that the spectrum $\s(\mL)$ of $\mL$
is absolutely continuous and consists of non-degenerated
intervals $[\l_{n-1}^+,\l_n^-], n=1,..,N_{G}\le \iy$.  
These intervals are separated by the gaps
$\g_n=(\l_n^-,\l_n^+)$ with the length $>0$.
{\bf Without loss of generality  we assume}
\[
\l_0^+=0,\ \ \ \ \ \ {\rm and} \ \ \ 
\lb{11} V^0=\int_0^1V(t)dt={\rm diag} \{V_{1}^0,
V_{2}^0,...,V_{N}^0\},\ \ \ \
V_{1}^0\le V_{2}^0\le ...\le V_{N}^0.
\]

A great number of  papers is devoted to the inverse spectral
theory for the Hill operator. We mention all papers where the
inverse problem including characterization was solved: Marchenko
and Ostrovski [MO1], Garnett and Trubowitz [GT1-2],Kappeler [Kap],
Kargaev and Korotyaev [KK1], and Korotyaev [K1-3] and for $2\ts 2$
Dirac operator Misura [Mi1-2] and Korotyaev [K4-5]. Recently, the author [K6]
extended the results of [MO1], [GT1], [K1-2] for the case $-y''+uy
$ to the case of distributions, i.e. $-y''+u'y$ on $L^2(\R)$,
where periodic $u\in L_{loc}^2(\R)$. It is important that in these
papers new results from analytic function theory (in
particular, conformal mapping theory) were obtained. As an 
example, we mention the proofs by the direct method (see [GT2],[KK1],[K1-3]). These
are short, but this approach needs a priori estimates of potential
in terms of spectral data. A priori estimates for various
parameters of the Hill operator
and for the Dirac operator (the norm of a periodic potential, effective masses,
gap lengths, height of slits, and so on) were obtained in [GT1],
[MO1-2], [KK1], [KK2], [K2-11].
 In order to get the required estimates  the authors of [GT1],
[MO1-2], [Mi1-2], [K2-11]... used the ``global quasi-momentum''(the conformal mapping), which was
introduced into the spectral  theory of  the  Hill operator by
Marchenko-Ostrovski [MO1].

There exist many papers about the periodic systems $N\ge 2$ (see [Ca1-3], [YS]). The basic results for direct spectral theory for the matrix case were obtained by Lyapunov [Ly]  (see also interesting papers of Krein [Kr], Gel'fand and Lidskii [GL]).
In [BBK] for the case $N=2$ the following results are obtained:
the Lyapunov function is constructed on the 2-sheeted 
Riemann surface and the existence of real and complex branch 
points is proved. In [BK] the operator $y''''+qy$, where
$q$ is a periodic real potential, was studied. In this case
the Lyapunov function is constructed on a 2-sheeted 
Riemann surface and the existence of real and complex branch 
points is proved. The asymptotics of gaps and resonances
in terms of the Fourier coefficients are obtained.

The main goal of our paper is to reformulate some spectral problem for the differential
operator with periodic matrix coefficients as problems of conformal mapping theory. We
construct the conformal mapping (averaged quasimomentum) $w$, with real part given by the
integrated density of states and imaginary part given by the Lyapunov exponent. We obtain
various properties of this conformal mapping, which are similar to the case $N=1$. For solving
these "new" problems we use some techniques from [KK2], [K2], [K6-8]. In particular, we use
the Poisson integral for the domain $\C_+\cup (-1,1)\cup\C_-$ and the Dirichlet integral for
the function $w(z)-z$. Note that the Dirichlet integral was used in [K6-8] for the scalar case to obtain a priori two-sided
estimates of the potential in terms of spectral data.

 Introduce the fundamental $N\ts N$-matrix solutions $\vp(t,z)$, $\vt(t,z)$ of the equation
\[
\lb{12}      -f''+V(t)f=z^2f,\ \ \ z\in\C,
\]
with the conditions $\vp(0,z)=\vt'(0,z)=0$,
$\vp'(0,z)=\vt(0,z)=I_N$, where $I_N,N\ge 1$ is the identity $N\ts
N$ matrix. Here and below we use the notation
$(')=\pa /\pa t$. We define the monodromy
$2N\ts 2N$-matrix $M$ and the trace $T_m, m\ge 1$ by
\[
\lb{13} M(z)=\mM(1,z),\ \ \ \mM(t,z)= \ma \vt(t,z)&\vp(t,z)\\
\vt'(t,z)&\vp'(t,z) \am, \ \  \ \ \ \ T_m(z)={\Tr M^m(z)\/2N}.
\]
The functions $M(z)$ and $T_m, m\ge 1$ are entire, real for $z^2\in \R$ and $\det M=1$. Let $\t_m, m=1,..,2N$ be the
eigenvalues of $M$. An eigenvalue of $M(z)$ is called a {\it multiplier}. It is a root of the algebraic equation 
$D(\t,z)\ev\det (M(z)-\t I_{2N})=0, \t,z\in\C$.
The zeros of $D(1,\sqrt\l)$ ( and $D(-1,\sqrt\l)$ (counted
with multiplicity) are the periodic (anti-periodic) eigenvalues
for the equation $-y''+Vy=\l y$ with periodic (anti-periodic) boundary conditions.

Below we need the following well-known results of Lyapunov (see [YS]).

{\bf Theorem (Lyapunov )} {\it Let $V\in \mH$. Then 
the following identities are fulfilled:
\[
\lb{TL-1}
 M^{-1}=-JM^\top J=\ma \vp'(1,\cdot)^\top & -\vp(1,\cdot)^\top\\
-\vt'(1,\cdot)^\top & \vt(1,\cdot)^\top\am,
\ \ \ \ \ J=\ma 0&I_N\\-I_N&0\am ,
\]
\[
\lb{TL-2}
 D(\t,\cdot)=\t^{2N}D(\t^{-1},\cdot),\ \ \ \t\neq 0.
\]
If for some $z\in\C$ (or $z^2\in\R$) $\t(z)$ is a multiplier
of multiplicity $d\ge 1$, then $\t^{-1}(z)$ (or $\ol\t(z)$) 
is a multiplier of multiplicity $d$. Moreover, each 
$M(z), z\in\C$, has exactly $2N$ multipliers
$\t_m^{\pm 1}(z), m=1,..,N$. Furthermore, $z^2\in \s(V)$ iff $|\t_m(z)|=1$  for some $m=1,..,N$.
If $\t(z)$ is a simple multiplier and $|\t(z)|=1$, then
$\t'(z)\ne 0$}.

It is well known that $ D(\t,z)=\sum_0^{2N}\x_m(z)\t^{2N-m}$,
where the functions $\x_m$ are given by
\[
\lb{15}
\x_0=1,\ \ \ \x_1=-2NT_1,\ \  \x_2=-{2N\/2}(T_2+T_1\x_1),\ \ ..
... ,
\x_m=-{2N\/m}\sum_0^{m-1}T_{m-j}\x_j,.. 
\]
see [RS]. Using the identity \er{TL-2} we obtain
\[
\lb{poD}
 D(\t,\cdot)=(\t^{2N}+1)+\x_1(\t^{2N-1}+\t)+
...+\x_{N-1}(\t^{N+1}+\t^{N-1})+\x_{N}\t^N.
\]
The eigenvalues of $M(z)$ are the zeros of Eq. $D(\t,z)=0$. This is an algebraic equation in $\t$ of degree $2N$. The coefficients $\x_m(z)$ are entire in $z\in\C$. It is well known
(see e.g. [Fo],[Sp]) that the roots $\t_m(z),m=1,..,2N$
constitute one or several branches of one or several analytic functions that have only algebraic singularities in $\C$.
Thus the number of eigenvalues of $M(z)$ is a constant $N_e$
with the exception of some special values of $z$
(see below the definition of a resonance). 
In general, there is a infinite number of such points
on the plane. If the functions $\t_m(z),m=1,..,N$ are all distinct, then $N_e=2N$. If some of them are identical,
then we get $N_e<2N$ and $M(z)$ is permanently degenerate. 

The Riemann surface for the multipliers $\t_m(z),m=1,..,N$  has $2N$ sheets, since degree of $D(\t,\cdot )$ is $2N$   see \er{poD}. If $N=1$, then it has 2 sheets, but the Lyapunov function is entire. Similarly, in the case $N\ge 2$ it is more convenient for us to  construct the Riemann surface for the Lyapunov function, which has N sheets (see Eq. \er{T1-2}). 
In order to  formulate our first result we 
transform $D(\t,z)$ to the polynomial $\F(\n,z)$ by
\[
\lb{14} {D(\t,z)\/(2\t)^N}=\F(\n,z)=\n^N+\f_1(z)\n^{N-1}+...+\f_N(z),\ \ \ \ \ \n={\t+\t^{-1}\/2},
\]
where $\f_1,..\f_N$ are some linear combinations of $\x_0,..\x_N$, see \er{110}-\er{113}. In
particular, all coefficients $\f_1(z),..\f_N(z)$ are entire functions. Each zero of  $\F(\n,z)$ is a Lyapunov
function 
$$
\D_m(z)={1\/2}(\t_m(z)+\t_m^{-1}(z)), \ m=1,..,N.
$$
{\bf Remark.} We note that this reduction from the polynomial
$D$ with $\deg D=2N$ to the polynomial
$\F$ with $\deg \F=N$ is crucial for our analysis. It is based
on \er{TL-2}, which is a consequence of $M$ being a symplectic
matrix and on the identity \er{CP} for the Chebyshev polynomials. 

We need the following preliminary results

\begin{theorem} \lb{T1}
Let $V\in \mH$. Then there exist analytic functions $\wt\D_s, s=1,..,N_0\le N$ on the
$N_s$-sheeted Riemann surface $\mR_s, N_s\ge 1$ having the following properties:

\no i) There exist disjoint subsets $\o_s\ss \{1,..,N\},
s=1,..,N_0, \bigcup \o_s=\{1,..,N\}$ such that all branches of
$\wt\D_s,s=1,2,..,N_0$ have the form
$\D_j(z)={1\/2}(\t_j(z)+\t_j^{-1}(z)), \ j\in \o_s$. Moreover,
 for any $z,\t\in \C$
the following identities are fulfilled:
\[
\lb{T1-2} {D(\t,z)\/(2\t)^{N}}=\prod_1^{N_0} \F_s(\n,z),\
\ \ \ \ \F_s(\n,z)=\prod_{j\in \o_s}(\n-\D_j(z)),\ \ 
\n ={\t+\t^{-1}\/2},\ \ \ \t\ne 0,
\]
where the functions $\F_s(\n,z)$ are entire with respect to $\n,z\in\C$ and $\F_s(\n,z)\in \R$ for all $\n,z\in \R$. Moreover, if $\D_i=\D_j$ for
some $i\in \o_k, j\in \o_s$, then $\F_k=\F_s$ and $\wt \D_k=\wt
\D_s$.

\no ii) (The monotonicity property). Let 
some branch $\D_m,m=1,..,N$ be real analytic on some interval
$Y=(\a,\b)\ss\R$ and $-1<\D_m(z)<1$ for any $z\in Y$.
Then $\D_m'(z)\ne 0$ for each $z\in Y$ .

\no iii) Each function $\r_s,s=1,..,N_0$ given by \er{T1-4} is entire and  real on the real line,
\[
\lb{T1-4}
\r=\prod_{1}^{N_0}\r_s,\ \ \ 
\r_s(\cdot)=\!\!\!\!\prod_{i<j, i,j\in \o_s}\!\!\!\! (\D_i(\cdot)-\D_j(\cdot))^2.
\]
\no iv) Each gap $\g_n=(\l_n^-,\l_n^+), n\ge 1$ is a bounded interval and $\l_n^\pm$ are either periodic (anti-periodic) eigenvalues or real branch points of $\D_m$ (for some $m=1,..,N$)
which are zeros of $\r$ (below we call such points resonances).
\end{theorem}

{\bf Remark.} 1)  If $N_0=N$, then $\mR_m=\C$ and each function $\D_m, m=1,...,N$ is entire.

2) We have the following asymptotics (see Sect. 3)
\[
\lb{1asD} \D_m(z)=\cos z+{\sin z\/2z}V_{m}^0+O\lt({e^{|\Im
z|}\/z^2}\rt), \ \ \  \ \ \ m=1,..,N,\  \ |z|\to \iy.
\]
Then firstly, $\r$ is not a polynomial since $\r$ is bounded on 
$\R$. Secondly, if $V_{i}^0\ne V_{j}^0, i\ne j$, then \er{1asD} implies
 $\D_i\ne \D_j$.

3) Let the surface $\mR=\cup_1^{N_0} \mR_s$ be a union of the disjoint  Riemann surfaces $\mR_s$ and let $\D=\{\wt\D_s, s=1,.,N_0\}$ be the corresponding analytic function on $\mR$.
Let $\f:\mR\to \C$ be the projection from the
surface $\mR$ into the complex plane. 
 {\bf We\  set} $\z\in \mR$ and $z\in \C$. 
 
  \vskip 0.25cm
{\bf Definition.} {\it The number $z_0$ is a {\bf resonance} of $\mL$, if $z_0$ is a zero of $\r$ given by \er{T1-4}.}

 \vskip 0.25cm

Define the real matrices $V^{sn}=\{V^{sn}_{jk}\}, V^{cn}=\{V^{cn}_{jk}\}$ by
\[
\hat V^{(n)}=\{\hat V^{(n)}_{jk}\}=\hat V^{cn}+i\hat V^{sn}=\int_0^1\!\!V(t)e^{i2\pi nt}dt.
\]
Denote by $\l_{m}^{n,\pm} , n\ge 0,m\in \{1,2,..,N\}$ the eigenvalues of the periodic and anti-periodic problem for the equation $-f''+Vf=z^2 f$. The periodic eigenvalues ($n$ is even) satisfy
\[
\lb{per}
0\le \underbrace{\l_{1}^{0,+}\le \l_{2}^{0+}\le ...\le\l_{N}^{0,+}}_{n=0}\le 
\underbrace{\l_{1}^{2,-}\le \l_{1}^{2,+}\le
...\le\l_{N}^{2,-}\le\l_{N}^{2,+}}_{n=2}\le \l_{1}^{4,-}\le \l_{1}^{4,+}\le \dots
\]
Recall $\l_0^+=0$. The anti-periodic eigenvalues ($n$ is odd) satisfy
\[
\lb{ape} 0\le \underbrace{\l_{1}^{1,-}\le \l_{1}^{1,+}\le
...\le\l_{N}^{1,-}\le\l_{N}^{1,+}}_{n=1}\le \underbrace{\l_{1}^{3,-}\le \l_{1}^{3,+}\le  ...
\le\l_{N}^{1,+}}_{n=3}\le \l_{1}^{5,-}\le \dots
\]
If $V=0$, then $\l_{m}^{n,\pm}=(\pi n)^2, m=1,..,N$.
Let $z_{m}^{n,\pm}=\sqrt{\l_{m}^{n,\pm}}>0$ and $z_{m}^{-n,\pm}=-z_{m}^{n,\mp}, n\ge 0, m\in
\{1,2,..,N\}$. The zeros of $D(1,z)$ and $D(-1,z)$ (counted
with multiplicity)  have the forms
$z_{m}^{2n,\pm}$ and $z_{m}^{2n+1,\pm}$ $n\in \Z$.
Let $|A|$ denote the operator norm of the matrix $A$.

\begin{theorem}   \lb{T2}
 Let $V\in \mH$.  Then the periodic and anti-periodic eigenvalues
have the following asymptotics:
\[
\lb{T2-1} \l_{m}^{n,\pm}=(\pi n)^2+\z_{m}^{n,\pm}+O(n^{-1}), \ \ \
m=1,..,N,\ n\to \iy,
\]
where $\z_{m}^{n,\pm}, m=1,2,..,N$ are the
eigenvalues of the matrix $\ma V^0+\hat V^{cn}&\!\!\hat V^{sn}\\ 
\hat V^{sn}\!\!& V^{0}-\hat V^{cn} \am$.

Assume that $V_j^0\ne V_{j'}^0$ for all $j\ne {j'}\in \o_s$ for some $s=1,..,N_0$. Then the function $\r_s$ has the zeros
$z_\a^{n\pm}, \a=(j,{j'}), j< j', j,j'\in \o_s, n\in \Z\sm \{0\}$, which are real at large $n$ and satisfy
\[
\lb{T2-2}
z_\a^{n\pm}=\pi n+{V^0_{i}+V^0_{j}\/4 \pi n}+O\lt({|\hat V^{(n)}|\/n}+{1\/n^2}\lt),\ \ \ \a=({j,j'})\ \ {\rm as}\ \  n\to \iy.
\]
Let in addition $V_{1}^0<..<V_{N}^0$. Then for each $s=1,..,N_0$
and for large $n\to \iy$ there exists a system of real intervals (gaps) $\g_n^\a=(\l_\a^{n-},\l_\a^{n+}), g_\a^n=(z_\a^{n-},z_\a^{n+}), 
$ such that 
$$z_{j,j'}^{n\pm}=z_{j',j}^{n\pm}>0,\ \ \ \  
\l_{j,j}^{n\pm}=\l_j^{n\pm},\ \ \ \ 
(z_\a^{n\pm})^2=\l_\a^{n\pm}>0,\ \ \a=(j,j'),\  j,j'\in \o_s,
$$ 
$$
\l_{j,j_1}^{n-}\le \l_{j,j_1}^{n+}<\l_{j,j_2}^{n-}\le \l_{j,j_2}^{n+}<...<
\l_{j,j_{N_s}}^{n-}\le \l_{j,j_{N_s}}^{n+}, \ N_s=|\o_s|
$$
\[
\lb{T2-3}
(-1)^n\D_j(z)>1, z\in g_{j,j}^n,\ \ and \ \  \ \ \ 
\ol \D_{j'}(z)=\D_j(z), z\in g_{j,j'}^n, \ \ {\rm if } \ \ j\ne j'
\]
i) The branch $\D_j$ is real and is analytic on the set $(\pi n-{\pi\/2},\pi n+{\pi\/2})\sm\cup_{p\ne j} g_{j,p}^n$ and is not real on $\cup_{p\ne j} g_{(j,p)}^n$.

\no ii) If $z_\a^{n-}\ne z_\a^{n+}$ for some $\a={j,j'},j\ne j'$, then  $z_\a^{n\pm}$ is a simple branch point, i.e., of square root type (resonance) for the functions $\D_j, \D_{j'}$. If $z_\a^{n-}=z_\a^{n+}$, then $\D_j, \D_{j'}$ are analytic at $z_\a^{n\pm}$. 

\no iii) The following asymptotics are fulfilled:
\[
\lb{T2-4} \l_\a^{n\pm}=(\pi n)^2+{V^0_{i}+V^0_{j}\/2}\pm
|\hat V^{(n)}_\a|+O\lt(|\hat V^{(n)}|+{1\/n}\lt),\ \ \ 
\a=({j,j'}).
\]
\end{theorem}

{\bf Remark.} 1) If  $N=1$, then
the asymptotics \er{T2-4} are well known [Ti]. 2) We describe the  surface $\mR$ in the case $V_1^0<..<V_N^0$ and $N_0=1$ for large $z$.  To "build" the surface $\mR$ for large $z$, we
take $N$ replicas of the z-plane  and call them sheets $\cR_1,..,\cR_N$. Each $\cR_j$ is cut along the real interval $g_\a^n, \a=({j,j'}),j,j'=1,..,N, j\ne j'$.
The cut on each sheet  two edges; we label each edge with $\sharp$ or $\flat$.   Then attach the $\sharp$ edge of the cut $g_\a^n$ on $\cR_j$ to the $\flat$ edge of the same cut on $\cR_{j'}, j\ne j'$, and attach the $\flat$ edge of the cut on $\cR_j$ to the $\sharp$ edge of the same cut on
$\cR_{j'}$. Thus, whenever we cross the cut, we pass from one sheet
to the other.  3) If $V_1^0<..<V_N^0$, then  all resonances  are real at high
energy. The existence of low energy complex resonances for
specific  potentials was established in [BKK].

\begin{figure}\lb{fig}
\centering \includegraphics[height=0.95\textheight]{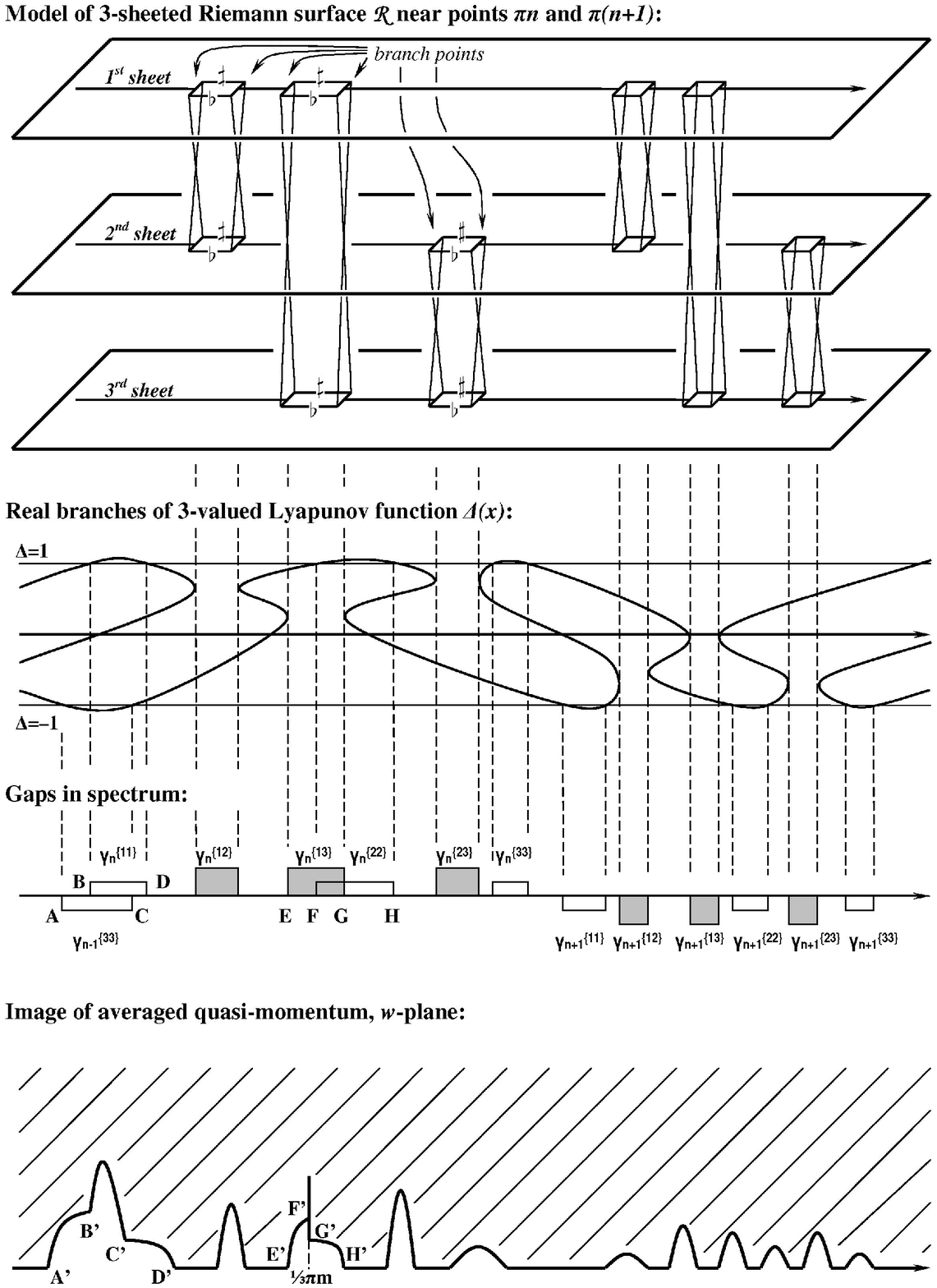} \caption{$N=3$,
$V_1^0<V_2^0<V_3^0$, $\frac{1}{2}\,(V_1^0+V_3^0)<V_2^0$.}
\end{figure}

Recall that $N_G$ is the total number of gaps in the spectrum of $\mL$.
\begin{corollary}\lb{Tc1} 
Let $V\in \mH$ and $V_{1}^0<...<V_{N}^0$.\\
(i) If the identity $V_1^0+V_N^0=V_2^0+V_{N-1}^0=...=V_N^0+V_1^0$ is not fulfilled, then $N_G<\iy$.\\
(ii) If $V_1^0+V_N^0=...=V_N^0+V_1^0$ holds true and there exists a sequence of indeces
$n_k\to\infty$ such that $|\hat V^{(n_k)}|^2\!+\!n_k^{-1}=o(|\hat V^{(n_k)}_{m,N+1-m}|)$ ,
$k\to\iy$, for each $m=1,..,N$, then $N_G=\iy$.
\end{corollary}
\noindent {\bf Remark.} 1) Maksudov and Veliev [MV] proved i),
for $N\ge 3$ in a more general case.\\
2) If $V=\diag\{V_{11},V_{22},..,V_{NN}\}$ and $V_1^0<V_2^0$, then the number of gaps is $N_G<+\infty$.\\
3) Note that the condition $|\hat V^{(n)}|^2\!+\!n^{-1}=o(|\hat V^{(n)}_{m,N+1-m}|)$,
$m=1,..,N$, $n\to\infty$, holds true for "generic"\ potentials from the space $\mH$. This yields the existence
of real resonance gaps $(\l_\a^{n-},\l_\a^{n+})$
at high energy. The coefficients $V_{m,N+1-m}, m=1,..,N$ 
(the "second diagonal" of the matrix $V$) "create" the gaps.

{\bf We consider the conformal mapping associated with the operator}
$\mL$. We need functions from the subharmonic  counterpart of the
Cartwright class of the entire functions given by
\[
\lb{SCdef} \cS\cC
 =\lt\{v: \C \to \R, \ \begin{array}{c} v {\rm \ is \ subharmonic \ in \ }
 \C {\rm \ and \ harmonic\ outside \ }\R, \cr
 v(\ol z)\equiv v(z),\ \int_{\R}{v_+(t)dt\/1+t^2}<\iy,\
 {\mathop{\lim\sup}\limits}_{z\to\iy}{v(z)\/|z|}<\iy\end{array}\rt\}.
\]
 We recall the class of functions from [KK1] given by
$$
\cS\cK_m^+=\lt\{v\in \cS\cC: v\ge 0, \ \ \lim_{y\to \iy}{v(iy)\/y}=1,\ \ \ 
\int_{\R}(1+t^{2m})v(t)dt<\iy\rt\}, \ \ m\ge 0.
$$
We note that $\cS\cK_{m+1}^+\ss\cS\cK_m^+, m\ge 0$.

Introduce the simple conformal mapping $\e:\C\sm [-1,1]\to \{\z\in
\C: |\z|>1\}$ by
\[
\lb{dx} \e(z)=z+\sqrt{z^2-1},\ \ \ \ z\in \C\sm [-1,1],\ \ {\rm
and}\ \ \ \  \e(z)=2z+o(1),\ \ |z|\to \iy.
\]
Note that $\e(z)=\ol \e(\ol z), z\in \C\sm [-1,1]$ since $\e(z)>1$
for any $z>1$. Due to the properties of the Lyapunov functions we have$|\e(\wt\D_s(\z))|>1, \z\in \mR_s^+=\{\z\in\mR_s: \Im \z>0\}$. Thus we can introduce the quasimomentum $k_m$ (we fix some branch of $\arccos$ and $\D_m(z)$) and the function $q_m$ by
\[
\lb{dkm} k_m(z)=\arccos \D_m(z)=i\log \e(\D_m(z)),\ \  
q_m(z)=\Im k_m(z)=\log |\e(\D_m(z))|,\ \ 
\]
$m=1,2,..,N$ and $z\in \mR_0^+=\C_+\sm \b_+,
\b_+=\!\bigcup_{\b\in \cB_\D\cap \C_+}
[\b, \b+i\iy)$
where $\cB_\D$ is the set of all branch points of the function $\D$.
The branch points of $k_m$ belong to $\cB_\D$. 
Define the {\bf averaged quasimomentum} $w$, the {\bf
density} $u$ and the {\bf Lyapunov exponent} $v$ by
\[
w(z)=u(z)+iv(z)={1\/N}\sum_1^N k_m(z), \ \ \  v(z)=\Im w(z),
\ \  z\in \mR_0^+.
\]
Define the sets $\s_{(N)}\!=\!\{z\in\R: \D_1(z),..,\D_N(z)\in [-1,1]\}$ and 
$$
\s_{(1)}\!=\!\{z\in \R: \D_m(z)\in (-1,1),\ \D_p(z)\notin [-1,1] \ {\rm \ some} \ m,p=1,..,N\}.
$$ 
For the function $w(z)=u(z)+iv(z), z=x+iy\in \ol\C_+$ we formally 
introduce the integrals
\[
\lb{dQSI} Q_n={1\/\pi}\int_{\R}x^nv(x)dx,\ \ \ \ \
 \cP_n={1\/\pi}\int_{\R}x^nv(x)du(x),\ \ \
 I_n^D={1\/\pi}\iint_{\C_+}|w_n'(z)|^2dxdy, \ \
\]
$n\ge 0$, where here and below $w_m(z),z\in\C_+$ is given by
$$
w_m(z)={1\/\pi}\int_{\R}{t^m v(t)\/t-z}dt =z^m\lt(w(z)-z+\sum_{n=0}^{m-1}Q_nz^{-n-1}\rt), \ \ \
z\in\C_+.
$$
 Let $C_{us}$ denote the class of all
real upper semi-continuous functions  $h:\R\to \R$. With any $h\in
C_{us}$  we associate the "upper" domain $W(h)=\{w=u+iv\in\C:
v>h(u), u\in \R\}$. Let $g=\cup_{n\in \Z} g_n$ where $g_n=(z_n^-,z_n^+), z_n^\pm=\sqrt{\l_n^\pm}>0$
and $g_{-n}=g_n, n\ge 1$. We formulate our main result. 

\begin{theorem}   \lb{T3}
Let $V\in \mH$. Then the averaged quasimomentum $w={1\/N}\sum_1^N k_m$ is analytic
in $\C_+$ and $w:\C_+\to w(\C_+)=W(h)$ is a conformal mapping
onto $W(h)$ for
some $h\in C_{us}$. Moreover, $v=\Im w$ has an harmonic extension
from $\C_+$ into $\O=\C_+\cup \C_-\cup g$ given by $v(z)=v(\ol z), z\in \C_-$ and $v(z)>0$ for any $z\in \O$. Furthermore  $v\in \cS\cK_2^+\cap C(\C)$ and there exist branches $k_m,m=1,..N$
such that the following asymptotics, identities and estimates are
fulfilled:
 \[
 \lb{T3-1}
w(z)=-\ol{w(-\ol z)},\ z\in\ol\C_+,\ \
\]
\[
 \lb{T3-2}
w(z)-z=-{Q_0\/z}-{Q_2+o(1)\/z^3},\ \ \  {\rm  as}\ y>r|x|, \ \  y\to\iy,\ \
\ \ {\rm for \ any} \ r>0,\ 
\]
\[
\lb{T3-3} Q_0=I_0^D+\cP_0=\int_0^1 {\Tr V(t)dt\/2N},
\]
\[
\lb{T3-4} Q_2=I_1^D+\cP_2=\int_0^1 {\Tr V^2(t)dt\/2^3N},
\]
\[
 \lb{T3-5}
v|_{\s_{(N)}}=0,\ \ \ 0<v|_{\s_{(1)}\cup g}\le \sqrt{2Q_0},\ \ \
\]
\end{theorem}

{\bf Remark.} 1) Craig and Simon  [CS] proved that
the Lyapunov exponent is subharmonic in $\C$ for 
the Schr\"odinger operator $-{d^2\/dx^2}+V$ for
a large class of scalar  potentials.

\begin{figure}[htb]
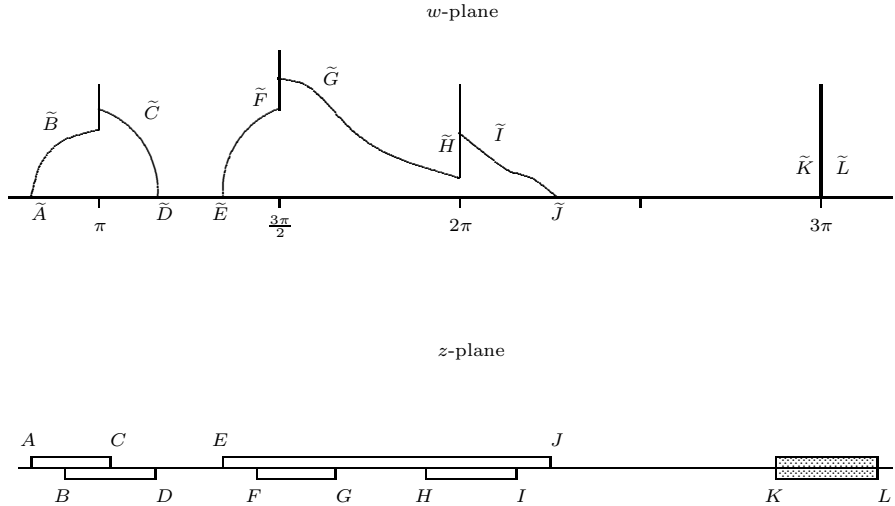

\begin{center}
\font\thinlinefont=cmr5
\begingroup\makeatletter\ifx\SetFigFont\undefined%
\gdef\SetFigFont#1#2#3#4#5{%
  \reset@font\fontsize{#1}{#2pt}%
  \fontfamily{#3}\fontseries{#4}\fontshape{#5}%
  \selectfont}%
\fi\endgroup%
\mbox{\beginpicture \setcoordinatesystem units <0.62992cm,0.62992cm> \unitlength=0.62992cm
\linethickness=1pt \setplotsymbol ({\makebox(0,0)[l]{\tencirc\symbol{'160}}}) \setshadesymbol
({\thinlinefont .}) \setlinear \linethickness= 0.500pt \setplotsymbol ({\thinlinefont .})
{\circulararc 72.861 degrees from  6.638 22.824 center at  7.300 21.129
}%
%
%
\linethickness= 0.500pt \setplotsymbol ({\thinlinefont .}) {\setshadegrid span <1pt>
\shaderectangleson \putrectangle corners at 17.145 15.479 and 19.287 15.240 \setshadegrid span
<5pt> \shaderectanglesoff
}%
%
%
\linethickness= 0.500pt \setplotsymbol ({\thinlinefont .}) {\setshadegrid span <1pt>
\shaderectangleson \putrectangle corners at 17.145 15.240 and 19.287 15.003 \setshadegrid span
<5pt> \shaderectanglesoff
}%
%
%
\linethickness= 0.500pt \setplotsymbol ({\thinlinefont .}) {\putrectangle corners at  5.476
15.479 and 12.383 15.240
}%
%
%
\linethickness= 0.500pt \setplotsymbol ({\thinlinefont .}) {\putrectangle corners at  9.762
15.240 and 11.667 15.003
}%
%
%
\linethickness= 0.500pt \setplotsymbol ({\thinlinefont .}) {\putrectangle corners at  6.191
15.240 and  7.857 15.003
}%
%
%
\linethickness= 0.500pt \setplotsymbol ({\thinlinefont .}) {\putrectangle corners at  2.142
15.240 and  4.047 15.003
}%
%
%
\linethickness= 0.500pt \setplotsymbol ({\thinlinefont .}) {\putrectangle corners at  1.429
15.479 and  3.095 15.240
}%
%
%
\linethickness= 0.500pt \setplotsymbol ({\thinlinefont .}) {\putrule from  1.156 15.240 to
19.6 15.240
}%
%
%
\linethickness= 0.500pt \setplotsymbol ({\thinlinefont .}) {\putrule from  0.953 20.955 to
19.6 20.955
}%
%
%
\linethickness= 0.500pt \setplotsymbol ({\thinlinefont .}) {\putrule from  6.638 23.467 to
6.644 23.467 \plot  6.644 23.467 6.659 23.465 / \plot  6.659 23.465  6.687 23.461 / \plot
6.687 23.461  6.723 23.457 / \plot  6.723 23.457  6.769 23.451 / \plot 6.769 23.451  6.822
23.442 / \plot  6.822 23.442  6.879 23.434 / \plot  6.879 23.434  6.936 23.423 / \plot  6.936
23.423  6.993 23.410 / \plot  6.993 23.410  7.046 23.398 / \plot  7.046 23.398 7.095 23.383 /
\plot  7.095 23.383  7.140 23.368 / \plot  7.140 23.368  7.184 23.351 / \plot  7.184 23.351
7.224 23.332 / \plot 7.224 23.332  7.262 23.309 / \plot  7.262 23.309  7.300 23.285 / \plot
7.300 23.285  7.338 23.258 / \plot  7.338 23.258  7.368 23.235 / \plot  7.368 23.235  7.400
23.207 / \plot  7.400 23.207 7.434 23.180 / \plot  7.434 23.180  7.468 23.148 / \plot  7.468
23.148  7.501 23.116 / \plot  7.501 23.116  7.537 23.080 / \plot 7.537 23.080  7.573 23.044 /
\plot  7.573 23.044  7.612 23.004 / \plot  7.612 23.004  7.650 22.964 / \plot  7.650 22.964
7.688 22.921 / \plot  7.688 22.921  7.728 22.879 / \plot  7.728 22.879 7.766 22.837 / \plot
7.766 22.837  7.806 22.792 / \plot  7.806 22.792  7.846 22.748 / \plot  7.846 22.748  7.885
22.705 / \plot 7.885 22.705  7.925 22.663 / \plot  7.925 22.663  7.963 22.623 / \plot  7.963
22.623  7.999 22.583 / \plot  7.999 22.583  8.037 22.543 / \plot  8.037 22.543  8.073 22.507 /
\plot  8.073 22.507 8.109 22.471 / \plot  8.109 22.471  8.145 22.435 / \plot  8.145 22.435
8.181 22.401 / \plot  8.181 22.401  8.217 22.367 / \plot 8.217 22.367  8.255 22.333 / \plot
8.255 22.333  8.293 22.301 / \plot  8.293 22.301  8.333 22.267 / \plot  8.333 22.267  8.376
22.236 / \plot  8.376 22.236  8.418 22.204 / \plot  8.418 22.204 8.460 22.172 / \plot  8.460
22.172  8.507 22.140 / \plot  8.507 22.140  8.551 22.109 / \plot  8.551 22.109  8.600 22.077 /
\plot 8.600 22.077  8.647 22.047 / \plot  8.647 22.047  8.695 22.020 / \plot  8.695 22.020
8.744 21.990 / \plot  8.744 21.990  8.791 21.965 / \plot  8.791 21.965  8.839 21.937 / \plot
8.839 21.937 8.888 21.912 / \plot  8.888 21.912  8.937 21.888 / \plot  8.937 21.888  8.983
21.865 / \plot  8.983 21.865  9.032 21.844 / \plot 9.032 21.844  9.078 21.823 / \plot  9.078
21.823  9.127 21.802 / \plot  9.127 21.802  9.172 21.785 / \plot  9.172 21.785  9.218 21.766 /
\plot  9.218 21.766  9.267 21.747 / \plot  9.267 21.747 9.318 21.728 / \plot  9.318 21.728
9.370 21.709 / \plot  9.370 21.709  9.426 21.689 / \plot  9.426 21.689  9.487 21.668 / \plot
9.487 21.668  9.553 21.647 / \plot  9.553 21.647  9.622 21.624 / \plot  9.622 21.624  9.696
21.601 / \plot  9.696 21.601  9.777 21.575 / \plot  9.777 21.575  9.859 21.550 / \plot  9.859
21.550 9.946 21.524 / \plot  9.946 21.524 10.031 21.497 / \plot 10.031 21.497 10.116 21.471 /
\plot 10.116 21.471 10.196 21.448 / \plot 10.196 21.448 10.270 21.425 / \plot 10.270 21.425
10.331 21.408 / \plot 10.331 21.408 10.382 21.391 / \plot 10.382 21.391 10.420 21.380 / \plot
10.420 21.380 10.446 21.374 / \plot 10.446 21.374 10.458 21.370 / \plot 10.458 21.370 10.465
21.368 /
}%
%
%
\linethickness= 0.500pt \setplotsymbol ({\thinlinefont .}) {\putrule from  6.668 22.80 to
6.668 24.052
}%
%
%
\linethickness=1pt \setplotsymbol ({\makebox(0,0)[l]{\tencirc\symbol{'160}}}) {\putrule from
18.098 20.955 to 18.098 23.336
}%
%
%
\linethickness= 0.500pt \setplotsymbol ({\thinlinefont .}) {\plot 10.465 22.316 10.471 22.310
/ \plot 10.471 22.310 10.486 22.299 / \plot 10.486 22.299 10.511 22.278 / \plot 10.511 22.278
10.547 22.250 / \plot 10.547 22.250 10.592 22.214 / \plot 10.592 22.214 10.645 22.172 / \plot
10.645 22.172 10.698 22.130 / \plot 10.698 22.130 10.753 22.085 / \plot 10.753 22.085 10.806
22.043 / \plot 10.806 22.043 10.854 22.005 / \plot 10.854 22.005 10.899 21.969 / \plot 10.899
21.969 10.941 21.937 / \plot 10.941 21.937 10.977 21.907 / \plot 10.977 21.907 11.013 21.880 /
\plot 11.013 21.880 11.045 21.855 / \plot 11.045 21.855 11.074 21.831 / \plot 11.074 21.831
11.102 21.808 / \plot 11.102 21.808 11.134 21.783 / \plot 11.134 21.783 11.165 21.759 / \plot
11.165 21.759 11.197 21.734 / \plot 11.197 21.734 11.229 21.711 / \plot 11.229 21.711 11.261
21.685 / \plot 11.261 21.685 11.292 21.662 / \plot 11.292 21.662 11.324 21.639 / \plot 11.324
21.639 11.356 21.618 / \plot 11.356 21.618 11.388 21.596 / \plot 11.388 21.596 11.417 21.577 /
\plot 11.417 21.577 11.447 21.558 / \plot 11.447 21.558 11.474 21.541 / \plot 11.474 21.541
11.502 21.526 / \plot 11.502 21.526 11.527 21.514 / \plot 11.527 21.514 11.553 21.503 / \plot
11.553 21.503 11.576 21.493 / \plot 11.576 21.493 11.604 21.482 / \plot 11.604 21.482 11.633
21.471 / \plot 11.633 21.471 11.663 21.463 / \plot 11.663 21.463 11.692 21.457 / \plot 11.692
21.457 11.722 21.448 / \plot 11.722 21.448 11.756 21.440 / \plot 11.756 21.440 11.788 21.433 /
\plot 11.788 21.433 11.822 21.425 / \plot 11.822 21.425 11.853 21.414 / \plot 11.853 21.414
11.887 21.406 / \plot 11.887 21.406 11.919 21.395 / \plot 11.919 21.395 11.949 21.385 / \plot
11.949 21.385 11.980 21.372 / \plot 11.980 21.372 12.012 21.357 / \plot 12.012 21.357 12.040
21.342 / \plot 12.040 21.342 12.067 21.325 / \plot 12.067 21.325 12.097 21.306 / \plot 12.097
21.306 12.131 21.283 / \plot 12.131 21.283 12.167 21.258 / \plot 12.167 21.258 12.207 21.226 /
\plot 12.207 21.226 12.251 21.192 / \plot 12.251 21.192 12.298 21.154 / \plot 12.298 21.154
12.349 21.114 / \plot 12.349 21.114 12.397 21.074 / \plot 12.397 21.074 12.442 21.035 / \plot
12.442 21.035 12.480 21.004 / \plot 12.480 21.004 12.507 20.983 / \plot 12.507 20.983 12.522
20.968 / \plot 12.522 20.968 12.531 20.961 /
}%
%
%
\linethickness= 0.500pt \setplotsymbol ({\thinlinefont .}) {\putrule from 10.478 21.368 to
10.478 23.336
}%
%
%
\linethickness= 0.500pt \setplotsymbol ({\thinlinefont .}) {\putrule from  2.857 22.384 to
2.855 22.384 \plot  2.855 22.384 2.849 22.382 / \plot  2.849 22.382  2.832 22.377 / \plot
2.832 22.377  2.800 22.371 / \plot  2.800 22.371  2.758 22.360 / \plot 2.758 22.360  2.703
22.348 / \plot  2.703 22.348  2.639 22.333 / \plot  2.639 22.333  2.574 22.316 / \plot  2.574
22.316  2.508 22.299 / \plot  2.508 22.299  2.445 22.280 / \plot  2.445 22.280 2.385 22.265 /
\plot  2.385 22.265  2.330 22.248 / \plot  2.330 22.248  2.282 22.231 / \plot  2.282 22.231
2.237 22.214 / \plot 2.237 22.214  2.197 22.200 / \plot  2.197 22.200  2.159 22.183 / \plot
2.159 22.183  2.125 22.166 / \plot  2.125 22.166  2.091 22.147 / \plot  2.091 22.147  2.057
22.128 / \plot  2.057 22.128 2.026 22.106 / \plot  2.026 22.106  1.994 22.083 / \plot  1.994
22.083  1.962 22.060 / \plot  1.962 22.060  1.930 22.035 / \plot 1.930 22.035  1.901 22.007 /
\plot  1.901 22.007  1.869 21.977 / \plot  1.869 21.977  1.839 21.946 / \plot  1.839 21.946
1.812 21.914 / \plot  1.812 21.914  1.782 21.882 / \plot  1.782 21.882 1.757 21.848 / \plot
1.757 21.848  1.731 21.816 / \plot  1.731 21.816  1.708 21.783 / \plot  1.708 21.783  1.687
21.749 / \plot 1.687 21.749  1.666 21.717 / \plot  1.666 21.717  1.649 21.685 / \plot  1.649
21.685  1.632 21.654 / \plot  1.632 21.654  1.615 21.620 / \plot  1.615 21.620  1.600 21.588 /
\plot  1.600 21.588 1.587 21.556 / \plot  1.587 21.556  1.573 21.520 / \plot  1.573 21.520
1.560 21.484 / \plot  1.560 21.484  1.547 21.444 / \plot 1.547 21.444  1.535 21.402 / \plot
1.535 21.402  1.522 21.355 / \plot  1.522 21.355  1.509 21.304 / \plot  1.509 21.304  1.494
21.249 / \plot  1.494 21.249  1.482 21.194 / \plot  1.482 21.194 1.469 21.137 / \plot  1.469
21.137  1.456 21.086 / \plot  1.456 21.086  1.446 21.040 / \plot  1.446 21.040  1.439 21.002 /
\plot 1.439 21.002  1.433 20.976 / \plot  1.433 20.976  1.431 20.961 / \plot  1.431 20.961
1.429 20.955 /
}%
%
%
\linethickness= 0.500pt \setplotsymbol ({\thinlinefont .}) {\putrule from  2.857 22.384 to
2.857 23.336
}%
%
%
\linethickness= 0.500pt \setplotsymbol ({\thinlinefont .}) {\putrule from 18.098 20.955 to
18.098 20.743
}%
%
%
\linethickness= 0.500pt \setplotsymbol ({\thinlinefont .}) {\putrule from 14.287 20.955 to
14.287 20.743
}%
%
%
\linethickness= 0.500pt \setplotsymbol ({\thinlinefont .}) {\putrule from 10.478 20.955 to
10.478 20.743
}%
%
%
\linethickness= 0.500pt \setplotsymbol ({\thinlinefont .}) {\putrule from  6.668 20.955 to
6.668 20.743
}%
%
%
\linethickness= 0.500pt \setplotsymbol ({\thinlinefont .}) {\putrule from  2.857 20.955 to
2.857 20.743
}%
%
%
\put{\SetFigFont{7}{8.4}{\rmdefault}{\mddefault}{\updefault}{$\wt L$}%
} [lB] at 18.404 21.465
%
%
\put{\SetFigFont{7}{8.4}{\rmdefault}{\mddefault}{\updefault}{$\wt K$}%
} [lB] at 17.554 21.465
%
%
%
%
\put{\SetFigFont{7}{8.4}{\rmdefault}{\mddefault}{\updefault}{$w$-plane}%
} [lB] at  9.762 24.765
%
%
\put{\SetFigFont{7}{8.4}{\rmdefault}{\mddefault}{\updefault}{$z$-plane}%
} [lB] at 10.001 17.621
%
%
%
%
%
%
%
%
%
%
%
%
\put{\SetFigFont{7}{8.4}{\rmdefault}{\mddefault}{\updefault}{$J$}%
} [lB] at 12.383 15.716
%
%
\put{\SetFigFont{7}{8.4}{\rmdefault}{\mddefault}{\updefault}{$I$}%
} [lB] at 11.667 14.527
%
%
\put{\SetFigFont{7}{8.4}{\rmdefault}{\mddefault}{\updefault}{$H$}%
} [lB] at  9.525 14.527
%
%
\put{\SetFigFont{7}{8.4}{\rmdefault}{\mddefault}{\updefault}{$G$}%
} [lB] at  7.857 14.527
%
%
\put{\SetFigFont{7}{8.4}{\rmdefault}{\mddefault}{\updefault}{$F$}%
} [lB] at  5.952 14.527
%
%
\put{\SetFigFont{7}{8.4}{\rmdefault}{\mddefault}{\updefault}{$D$}%
} [lB] at  4.047 14.527
%
%
\put{\SetFigFont{7}{8.4}{\rmdefault}{\mddefault}{\updefault}{$B$}%
} [lB] at  1.905 14.527
%
%
\put{\SetFigFont{7}{8.4}{\rmdefault}{\mddefault}{\updefault}{$A$}%
} [lB] at  1.190 15.716
%
%
\put{\SetFigFont{7}{8.4}{\rmdefault}{\mddefault}{\updefault}{$C$}%
} [lB] at  3.095 15.716
%
%
\put{\SetFigFont{7}{8.4}{\rmdefault}{\mddefault}{\updefault}{$E$}%
} [lB] at  5.239 15.716
%
%
\put{\SetFigFont{7}{8.4}{\rmdefault}{\mddefault}{\updefault}{$L$}%
} [lB] at 19.287 14.527
%
%
\put{\SetFigFont{7}{8.4}{\rmdefault}{\mddefault}{\updefault}{$K$}%
} [lB] at 16.906 14.527
%
%
\put{\SetFigFont{7}{8.4}{\rmdefault}{\mddefault}{\updefault}{$\wt J$}%
} [lB] at 12.383 20.479
%
%
\put{\SetFigFont{7}{8.4}{\rmdefault}{\mddefault}{\updefault}{$\wt I$}%
} [lB] at 11.191 22.147
%
%
\put{\SetFigFont{7}{8.4}{\rmdefault}{\mddefault}{\updefault}{$\wt H$}%
} [lB] at 10.001 21.907
%
%
\put{\SetFigFont{7}{8.4}{\rmdefault}{\mddefault}{\updefault}{$\wt G$}%
} [lB] at  7.586 23.338
%
%
\put{\SetFigFont{7}{8.4}{\rmdefault}{\mddefault}{\updefault}{$\wt F$}%
} [lB] at  6.088 22.9
\put{\SetFigFont{7}{8.4}{\rmdefault}{\mddefault}{\updefault}{$\wt E$}%
} [lB] at  5.239 20.479
%
%
\put{\SetFigFont{7}{8.4}{\rmdefault}{\mddefault}{\updefault}{$\wt C$}%
} [lB] at  3.810 22.623
%
%
\put{\SetFigFont{7}{8.4}{\rmdefault}{\mddefault}{\updefault}{$\wt B$}%
} [lB] at  1.666 22.384
%
%
\put{\SetFigFont{7}{8.4}{\rmdefault}{\mddefault}{\updefault}{$\wt D$}%
} [lB] at  4.047 20.479
%
%
\put{\SetFigFont{7}{8.4}{\rmdefault}{\mddefault}{\updefault}{$\wt A$}%
} [lB] at  1.429 20.479
%
%
\put{\SetFigFont{7}{8.4}{\rmdefault}{\mddefault}{\updefault}{$\pi$}%
} [lB] at  2.75 20.253
%
%
\put{\SetFigFont{7}{8.4}{\rmdefault}{\mddefault}{\updefault}
{$\frac{3\pi}{2}$}%
} [lB] at  6.25 20.253
%
%
\put{\SetFigFont{7}{8.4}{\rmdefault}{\mddefault}{\updefault}{$2\pi$}%
} [lB] at 10.25 20.253
%
%
\put{\SetFigFont{7}{8.4}{\rmdefault}{\mddefault}{\updefault}{$3\pi$}%
} [lB] at 17.858 20.253
%
\linethickness= 0.500pt \setplotsymbol ({\thinlinefont .}) {\circulararc 76.531 degrees from
4.098 20.961 center at  2.290 21.099
}%
\linethickness=0pt \putrectangle corners at  0.906 25.146 and 21.002 14.385
\endpicture}
\end{center}

\caption{$N=2$. The domain $W(h)=w(\C_+)$ and gaps in the spectrum.
} \label{fig1}
\end{figure}

\no  2) Similar arguments give
\[
\lb{T3-6} Q_4=I_2^D+\cP_4=\int_0^1 {\Tr (V'(t)^2+2V^3(t))dt\/2^5N},
\ \ \ \ \ if \ \ V'\in \mH.
\]
3) The integral $\cP_0\ge 0$ is the area between the boundary of $W(h)$  and the real line.
The mapping $w:\C_+\to W(h)$ is illustrated in Figure \ref{fig1}.
 In Figure 1 the upper picture is a domain
$W(h)$ and $\wt A=w(A), \wt B=w(B),...$  . The
spectral interval $(A,B)$ (with multiplicity 2) of the z-domain
is mapped on the curve $(\wt A,\wt B)$ of the w-domain, the
interval (a gap)  $(B,C)$ of the z-domain  is mapped on a
vertical slit, which lies on the line $\Re w=\pi$. The spectral
interval $(C,D)$ (with multiplicity 2) of the z-domain is mapped
on the curve $(\wt C,\wt D)$ of the w-domain. The spectral
interval $(D,E)$ (with multiplicity 4) of the z-domain is mapped
on  the interval $(\wt D,\wt C)$ of the w-domain. The case of the
interval $(E,J)$ is similar. The resonance gap $(K,L)$  of the
z-domain is mapped on the vertical slit on the line
$\Re w=3\pi$. In fact the boundary of $W(h)$
is given by the graph of the function $h(u),u\in \R$.


Let $\s(m,V)$ denote the spectrum of $\mL$ of multiplicity $2m,m\ge 0$.  We have the following simple corollary from Theorem  \ref{T3}.

\begin{corollary}   \lb{T4}
{\it Let $\s(\mL)=\s(N,V)=\R_+$ for some $V\in \mH$. Then
$V=0$.}
\end{corollary}

\no {\it Proof.} Due to $\s_{(N)}=\R_+,$ we obtain $Q_0=Q_2=0$. Then
identity  \er{T3-3} yields $\|V\|=0$.
 $\BBox$

Recall that in the scalar case,  the so-called Borg Theorem follows
immediately from the existence of the conformal mapping and the
asymptotics of the Lyapunov function at high energy [M]. In general,
in order to prove the uniqueness result (the simplest part in the
inverse spectral theory) it is necessary to use some results from
 "function theory". In our case we use conformal mapping
theory. Recall that the so-called Borg Theorem for periodic systems
was proved in [CHGL],[GKM] for general cases.

We describe the properties of the conformal mapping $w$.

\begin{theorem}   \lb{T5}
Let $V\in \mH$. Then the following relations are fulfilled:
 \[
 \lb{T5-1}
u_x'(z)\ge 1,\ \ \  z\in\s_{(N)}\ \ {\rm and}\ \ u_x'(z)>0,\ \ \
\ \ z\in\s_{(1)},
\]
here $u_x'(z)=1$ for some $z\in\s_{(N)}$ iff $V=0$. Moreover,
\[
\lb{T5-2}
 v_{xx}''(z)<0<v(z),\ \ \ \   u(z)={\rm const}\in {\pi\/N}\Z,
\  {\rm for \ all} \ \ z\in g_n=(z_n^-,z_n^+),
\]
\[
\lb{T5-3} v(x)=v_n^0(x)\rt(1+{1\/\pi}\int_{\R\sm g_n} {v(t)dt\/
v_n^0(t)|t-x|}\rt), \ \ x\in g_n,\ 
v_n^0(z)=|(z-z_n^-)(z_n^+-z)|^{1\/2},
\]
\[
 \lb{T5-4}
\sum_n |g_n|^2\le 8Q_0, \ \ \ \ \ \ \ \ \ \ 
G^2\ev\sum |\g_n|^2\le {8\/N}\|V\|^2,
\]
\[
\lb{T5-5} 
\|V\|\le C_0G(1+G^{1\/3}), \ \ \ {\rm if}\ \ \s(N,V)=\s(\mL),
\]
for some absolute constant $C_0$.
\end{theorem}

{\bf Remark.} Using this theorem we deduce that the function
$h(u)=v(x(u)), u\in \R$ is continuous on $\R\sm \{u_n, n\in \Z\}$,
where $u_n=u(x), x\in g_n$. In this case we have
\[
h(u_n\pm 0)\le h(u_n), \ \ \ n  \in \Z.
\]

The plan of our paper is as follows. In Sect. 2 we obtain the basic properties
of the fundamental solution and prove Theorem \ref{T1}. In Sect. 3 we determine the asymptotics  of
$M(z)$ and of the Lyapunov function and the multipliers at high energy and prove Theorem \ref{T2}. Sect. 4 is the central part of the paper where we obtain the main properties of the quasimomentum
$k_m,m=1,..,N$. We prove the basic Theorems \ref{T41} and
\ref{T42}. In Sect. 5
using Theorems \ref{T41}, \ref{T42} and [KK3], devoted to
the conformal mapping theory, we prove Theorem \ref{T3} and \ref{T5}.

\section {Fundamental solutions}
\setcounter{equation}{0}

In this section we study $\vt,\vp$. We begin with some notational
convention. A vector $h=\{h_n\}_1^N\in \C^N$ has the Euclidean
norm $|h|^2=\sum_1^N|h_n|^2$, while a $N\ts N$ matrix $A$ has the
operator norm given by $|A|=\sup_{|h|=1} |Ah|$.

The fundamental solution $\vp$ satisfy the integral equations
\[
 \lb{21}
\vp(t,z)= \vp_0(t,z)+\int_0^t{\sin z(t-s)\/z}V(s)\vp(s,z)ds,\ \ \
\vp_0(t,z)={\sin z t\/z}I_N,\ \ \
\]
where $(t,z)\in\R\ts\C$. The standard iterations in \er{21} yields
\[
\lb{22}
 \vp(t,z)={\sum}_{n\ge 0} \vp_n(t,z)\,, \quad
\vp_{n+1}(t,z)= \int_0^t{\sin z(t-s)\/z}V(s)\vp_n(s,z)ds.
\]
The similar expansion $\vt={\sum}_{n\ge 0} \vt_n$ with
$\vt_0(t,z)=I_N\cos zt$ holds.  In order to determine the
asymptotics of the fundamental solutions we introduce the
functions
\[
\lb{23} I_m^0(z)=\int_0^mds_1\int_0^{s_1}\cos
z(m-2s_1+2s_2)F_2(s)ds_2,\ \ \ \ \ \ m=1,2,
\]
where $F_n(s)=\Tr V(s_1)\cdot ..\cdot V(s_n),
s=(s_1,..,s_n)\in \R^n, n\ge 1$. They satisfy the simple identity
$I_2^0(z)=4I_1^0(z)\cos z,\ \ z\in \C$,
see [BKK]. Define
\[
\lb{25} V^{0}=\!\!\int_0^1\!\!\!\!V(x)dx,\ \ B_{n}=\Tr
{(V^{0})^n\/n!},\ \ n\ge 1,\ \ \    \ |z|_1=\max\{1,|z|\}.
\]
We prove

\begin{lemma} \lb{T21}
Let $V\in\mH$. Then each functions
$\vp(t,z),\vt(t,z),t\ge 0$, are entire and for any integers $m\ge t\ge 0, n_0\ge -1$ the following estimates are fulfilled:
$$
\max\lt\{\lt|\vt(t,z)-\sum_0^{n_0}\vt_n(t,z)\rt|,
\lt||z|_1\lt(\vp(t,z)-\sum_0^{n_0}\vp_n(t,z)\rt)\rt|,
\lt|{1\/|z|_1}\lt(\vt'(t,z)-\sum_0^{n_0}\vt_n'(t,z)\rt)\rt|,
$$
\[
\lb{27}
\lt|\vp'(t,z)-\sum_0^{n_0}\vp_n'(t,z)\rt|\rt\} \le
{(m\vk)^{n_0+1}\/(n_0+1)!}e^{t|\Im z|+m\vk}.
\]
Moreover, each $T_m(\cdot), m\ge 1$ is entire and  satisfies
\[
\lb{29}
 |T_m(z)-\sum_0^{n_0}T_m^{(n)}(z)|\le
 {(m\vk)^{n_0+1}\/(n_0+1)!}e^{m(|\Im z|+\vk)},
\]
$$ 
T_m^{(0)}(z)=\cos mz,\ \ \ \ T_m^{(1)}(z)= {mB_1\/2Nz}\,{\sin mz},\ \ \ \
T_m^{(2)}(z)={1\/4Nz^2}\,(I_2^{(m)}(z)-m^2B_2\cos mz),...
$$ 
\end{lemma}

\no {\it Proof.} We prove the estimates of $\vp$, the proof for
$\vp',\vt,\vt'$ is similar. \er{22} gives
\[
\lb{211}
\vp_n(t,z)= \int_{D_n}f_n(t,s)V(s_1)\cdot ..\cdot V(s_n)ds, \ \ \ \
\ f_n(t,s)=\vp_0(s_n,z)\prod_1^n {\sin z(s_{k-1}-s_{k})\/z},
\]
where $s=(s_1,..,s_n)\in \R^n, s_0=t$ and $D_n=\{0<s_n<
...<s_{2}<s_{1}<t\}$. 
Substituting the
estimate $|\vp_0(t,z)|=|z^{-1}\sin zt|\le |z|_1^{-1}e^{|\Im z|t}$ into \er{211}, we obtain
$$
|\vp_n(t,z)|\le{e^{|\Im z|t}\/|z|_1^{n+1}}\int_{D_n}|V(s_1)|\cdot ..\cdot |V(s_n)|ds\le
{e^{|\Im z|t}\/|z|_1^{n+1}}\cdot{1\/n!}\lt(\int_0^t|V(x)|dx\rt)^n.
$$
This shows that for each $t\ge 0$
the series \er{22} converges uniformly on bounded subset of $\C$.
Each term of this series is an entire function. Hence the sum is
an entire function. Summing the majorants we obtain estimates
 \er{27}. 
 
 The function $T_m,m \ge 1$ is entire, since $\vp,\vt$
 are entire. We have 
 $(2N)T_m=\Tr M^m(z)=\Tr \mM(m,z)=\Tr \sum_{n\ge
0}\mM_n(m,z)$, where
\[
\lb{212} \Tr \mM_0(m,z)=2N\cos mz,\ \ \ \ \Tr \mM_n(m,z)=\Tr
\vt_n(m,z)+\Tr \vp_n'(m,z),\ \ n\ge 1.
\]
The estimates $|\Tr \vp_n'(m,z)|\le {(m\vk)^n\/n!}e^{|\Im z|m}$ and
$|\Tr \vt_n(m,z)|\le {(m\vk)^n\/n!}e^{|\Im z|m}$ yield
\[
\lb{213}
|\Tr \mM_n(m,z)|\le (2N){(m\vk)^n\/n!}e^{m|\Im z|},\ \ n\ge 0.
\]
Using \er{212} we obtain
\[
\lb{214}
\Tr \mM_1(m,z)={1\/z}\int_0^m\!\!\! (\sin z(m-t)\cos zt+\cos
z(m-t)\sin zt )\Tr V(t)dt={\sin mz\/z}mB_1,
\]
 and
$$
\Tr \mM_2(m,z)={1\/z^2}\int_0^m\int_0^t\!\!\! \sin
z(t-s)z(m-t+s)F_2(t,s)dtds
$$
$$
={1\/2z^2}\int_0^m\int_0^t\!\!\! (\cos z(m-t+s)-\cos
zm)F_2(t,s)dtds= {1\/2z^2} (I_m^0(z)-m^2B_2\cos mz )
$$
since $
\int_0^m\int_0^t F_2(t,s)dtds={1\/2}\Tr
\lt(\int_0^mV(t)dt\rt)^2=m^2B_2.
$
 $\BBox$

Note that the fundamental solutions $\vp(t,z),\vt(x,z)$ and $M(z), T_m(z), m\ge 1$ is real for $z^2\in \R$.
Moreover, all functions $\vp(t,z),\vp'(t,z),\vt(t,z),\vt'(t,z), M(z),  T_m, m\ge 1$ are even with respect to $z\in \C$ and then they are entire with respect to $\l=z^2$. 

Using the identity \er{TL-2} and $ D(\t,z)=\sum_0^{2N}\x_m(z)\t^{2N-m}$ we obtain 
$$
{D(\t,z)\/(2\t)^N}=
{(\t^N+\t^{-N})\/2^N}+\x_1{(\t^{N-1}+\t^{1-N})\/2^N}+..+
\x_{N-1}{(\t+\t^{-1})\/2^N}+...+2^{-N}\x_{N}.
$$
Substituting into this equality the identity for the Chebyshev
polynomials $\cT_n, n\ge 1$:
\[
\lb{CP}
{\t^{n}+\t^{-n}\/2}=\cT_n(\n)=
2^{N-1}\sum_0^{[{n\/2}]}c_{n,m}\n^{n-2m},\ \ \ 
c_{n,m}=(-1)^mn{(n-m-1)!\/(n-2m)!m!}2^{n-2m-N},\ \ 
\]
$\n={\t+\t^{-1}\/2}$, see [AS], we get $\F(\n,z)$ given by
\[
\lb{110}
\F(\n,z)={D(\t,z)\/(2\t)^N}=\sum_0^N \f_m(z)\n^{N-m},\ \ \ \  
\ \ 
\]
\[
\lb{111}
\f_0=1,\ \ \f_1=c_{N-1,0}\x_1={\x_1\/2},\ \ \f_2=
c_{N,1}+c_{N-2,0}\x_2,\ \ \f_3=c_{N-1,1}\x_1+c_{N-3,0}\x_3,...,
\]
\[
\lb{112}
\f_{2n}=c_{N,n}+c_{N-2,n-1}\x_2+c_{N-4,n-2}\x_4+...+
c_{N-2n,0}\x_{2n},
\]
\[
\lb{113}
\f_{2n+1}=c_{N-1,n}\x_1+c_{N-3,n-1}\x_3+c_{N-5,n-2}\x_5+...+
c_{N-2n-1,0}\x_{2n+1}.
\]
Let $D^0, \F^0$ be the determinant $D$ and the polynomial $\F$ at $V=0$. In this case we have
$$
D^0(\t,z)=(\t^2+1-2\t\cos z)^N=(2\t)^N(\n-\cos z)^N
=(2\t)^N\sum_0^NC_m^N(-\cos z)^m\n^{N-m}
$$
where $C_m^{N}={N!\/(N-m)!m!}$.
Thus $\F^0=(\n-\cos z)^N$ and $\f_m^0(z)=(-1)^mC_m^N\cos^m z$ at $V=0$.

Substituting the identity for the Chebyshev
polynomials \er{CP} into the trace formula
${1\/2}\Tr M^n(z)=\sum_1^N{\t_m^n+\t_m^{-n}\/2}$,
we obtain 
\[
\lb{T1-1}
\sum_{m=1}^N\cT_n(\D_m(z))={1\/2}\Tr M^n(z),\ \ \ z\in \C.
\]

\no {\bf Proof of Theorem  \ref{T1}} 
$\D_1(z),..,\D_N(z)$ are the roots of $\F(\n,z)=0$ for fixed
$z\in \C$. Recall that $\cB_\D$ are all branch points of $\D_1(z),..,\D_N(z)$.

First simple case, let $\cB_\D=\es$. Then all functions $\D_1(z),..,\D_N(z)$ are entire, the function $\F_j=\n-\D_j, j=1,..,N$ and $N_0=N$.

Second case, let $\cB_\D\ne \es$.
Consider a simply-connected domain $\O_1\ss \C$ containing only
one branch point $z_1\in \cB_\D$. We consider the behavior of the roots $\D_m,m=1..,N$ in the neighborhood of the branch points $z_1$.
Let $B'(z_1,r)=B(z_1,r)\sm \{z_1\}$ be the small disk near $z_1$
with radius $r>0$  but excluding $z_1$.
 The functions $\D_m,m=1..,N$ are  
 branchs of analytic functions (defined in $\O_2\sm (\{z_1,z_2\})$)
 with a branch point (if $p_2>1$) at $z_1$ and $z_2$.

 are analytic of $z\in \cB'(z_1,r)$.
If $B'(z_1,r)$ is moved continuously around $z_1$, then $N$
functions can be continued analytically. When $B'(z_1,r)$ has
been brought to its initial position after one revolution around
$z_1$, the functions $\D_m,m=1..,N$ will have undergone a
permutation among themselves. These functions may therefore be
grouped in the manner
\[
\{\D_1(z),..,\D_{p_1}(z)\}, \{\D_{p_1+1}(z),..,\D_{p_1+q_1(z)}\},..
\ \  {\rm where }\ \ 
\D_i\ne \D_j, 1\le i<j\le p_1,...
\]
in such a way that each group undergoes a ciclic permutation by
a revolution of $B'(z_1,r_1)$ of the kind described.
For brevity each group will be called a cycle at the branch point
$z_1$, and the number of elements of a cycle will be called its period.

It is obvious that the elements of a cycle of period $p_1$ constitute a branch of an analytic function (defined near $z_1$) with a branch point (if $p_1>1$) at $z_1$. We have Puiseux series as
$$
\D_j(z)=\D_1(\z)+a_1t+a_2t^2+..,\ \ t=e^{j{i2\pi\/p_1}}(z-z_1)^{1\/p_1},\
j=1,..,p_1.
$$
Consider a simply-connected domain $\O_2\ss \C, $ containing only
 two distinct branch point $z_1, z_2\in \cB_\D$. The similar argument as above give that functions $\D_1,..,\D_N$ may therefore be grouped in the manner
\[
\{\D_1(z),..,\D_{p_2}(z)\}, \{\D_{p_2+1}(z),..,\D_{p_2+q_2(z)}\},
.. , {\rm where }\ \ p_1 \le p_2,..
\]
and $\D_i\ne \D_j, 1\le i<j\le p_1,...$.
Here the elements of a cycle of period $p_2$ constitute
a branch of an analytic function (defined in $\O_2\sm (\{z_1,z_2\})$)
 with a branch point (if $p_2>1$) at $z_1$ and $z_2$.

If we take a sequence of domain $\O_n\ss \O_{n+1}, n\ge 1$,
containing only $n$ distinct branch points $z_1,.., z_n\in \cB_\D\cap \O_n$ and let $ \lim \O_n=\C$, then the similar argument as above
give that functions $\D_1,..,\D_N$ may therefore be grouped in the manner
\[
\{\D_j(z), j\in \o_1\}, \{\D_j(z), j\in \o_2\},...,
\{\D_j(z), j\in \o_{N_0}\},\  \ {\rm where }\ \ 
\D_i\ne \D_j, i,j\in \o_s,
\]
$i\ne j, s=1,..,N_0$. 
Here the elements of $\{\D_j(z), j\in \o_1\}$ constitute
a branch of an analytic function $\wt \D_1$ (defined in $\C\sm \cB_\D$)
 with a branch point (if the of numbers of elements of $\o_1$ is $>1$) at some points $z_n \in \cB_\D, z\ge 1$.

There exists an interval $Y\ss \R, Y\neq \es$ such that
the spectral interval $Y\ss \s(\mL)$ has multiplicity $2N$
(see [MV]). Thus all functions $\D_m,m=1,..,N$ are real
on $Y$. Hence each entire function $\F_s(\n,z), s=1,..,N_0$ is real for
all $\n,z\in \R$. 

ii) We have $\D_m'(z)={1\/2}(1-\t^{-2}(z))\t'(z)\ne 0, z\in Y$, since by the Lyapunov Theorem, $\t'(z)\ne 0$ for all $z\in Y$.

iii) Recall  that the resultant for the polynomials 
$
f=\t^n+a_{1}\t^{n-1}+..+a_n, \  g=b_0\t^s+b_{1}\t^{s-1}+..+b_s
$
is given by
\[
\lb{res} 
R(f,g)=\det 
\left(\begin{array}{cccccccc}
1& a_1&..&a_{n}&0&0&..&0\\
0&1& a_1&..&a_{n}&0&..&0\\
. &.&.&.&. &.&.&. \\
0&...&0&1& a_1&&..&a_{n}\\
b_0&b_1&..&b_s&0&..&0\\
0&b_0&b_1&..&b_{s}&0&..&0\\
. &.&.&.&. &.&.&. \\
0&...&0&b_0&b_1&&..&b_s\\
\end{array}\right)
\begin{array}{c}
\left.\phantom{\begin{array}{c}
~\\~\\~\\~\end{array}}\right\} \mbox{ $s$ lines}\\
\left.\phantom{\begin{array}{c}
~\\~\\~\\~\end{array}}\right\}\mbox{ $n$ lines}
\end{array}.
\]
The discriminant of the polynomial $f$ with zeros $\t_1,..,\t_n$
is given by 
$${\rm Dis} f=\prod_{i<j} (\t_i-\t_j)^2=(-1)^{{n(n-1)\/2}}R(f,f').
$$  
Thus we have $\Dis \F_j(\t,z)=\!\!\prod_{i<s, i,s\in \o_j}\!\! (\D_i(z)-\D_s(z))^2=(-1)^{N_j(N_j-1)\/2}R(\F_j,\F_j')
$
is entire, since $\F_j(\t,z)$ is the entire function.
Then the function $\r=\prod_{1}^{N_0}\Dis \F_j$ is entire.

iv) We consider the gap $g_n=(z_n^-,z_n^+)$ in the variable $z$, where a gap $\g_n=(\l_n^-,\l_n^+), \l_n^\pm=(z_n^\pm)^2$.
Each gap $g_n=(z_n^-,z_n^+)=\cup g_{n,p}, p=1,..,p_n$, where
$g_{n,p}=(z_{n,p}^-,z_{n,p}^+)$ is finite interval such that
$\D_m(z)\notin [-1,1]$ for  all $z\in g_{n,p}$ for some $m=m(p)$.
Note that $\D_m(z_{n,p}^-)=\pm 1$ or $z_{n,p}^-$ is the branch
point $\in \cB_\D$, otherwise we have a contradiction.
$\BBox$

\section {Spectral asymptotics}
\setcounter{equation}{0}

Below we need the identities  for
$J=\ma 0&I_N\\-I_N&0\am , J_1=
\ma I_N&0\\0&-I_N\am , \ J_2= \ma 0&I_N\\I_N&0\am $:
\[
\lb{31}
  J^2=-I, \ \ \ J_1J_2=J,\ \ \ JJ_1=-J_2,\ \ \ JJ_2=J_1,
\]
\[
\lb{32} e^{zJ}=I_{2N}\cos z +J\sin z,         \ \ \ z\in \C.
\]
We have also
\[
\lb{33} D^0(z,\t)=\det (e^{zJ}-\t
I_{2N})=(\t-e^{iz})^N(\t-e^{-iz})^N,\ \ \ \ \
D^0(z,\pm 1)=2^N(1\mp\cos z)^N,
\]
\[
\lb{34} (e^{zJ}-\t I_{2N})^{-1}=
(\t-e^{iz})^{-1}(\t-e^{-iz})^{-1}(e^{-zJ}-\t I_{2N}).
\]
We shall obtain the simple properties of the monodromy matrix.
We introduce the modified monodromy matrix
\[
\lb{35} \wt M(z)=\ma I_N& 0\\ 0& zI_N\am^{-1}M(z) \ma I_N& 0\\ 0&
zI_N\am=\ma \vt(1,z)& z\vp(1,z)\\ z^{-1}\vt'(1,z)& \vp'(1,z)\am,
 \ \ \ z\in \C,
\]
with the same eigenvalues and the same traces.
We will show following asymptotics
\[
\lb{36} \wt M(z)=e^{zJ}\rt(I+{-V_0J+\hat V(z)J_2\/2z}\rt)+
O(z^{-2}e^{|\Im z|}),\ \ \hat V(z)=\!\!\int_0^1\!\!V(t)e^{-2tzJ}dt,
 |z|\to \iy.
\]
Indeed, using \er{27}, \er{35}  we get
$\wt M(z)=e^{zJ}+\wt M_1(z)+O(z^{-2}e^{|\Im z|})$, where $\wt M_1$ is given by
\[
\lb{37} \!\!\wt M_1(z)={1\/2z}\int_0^1\!\!\!\!V(t)f(t,z)dt,\
f=2\ma \sin z(1-t)\cos tz I_N& \sin z(1-t)\sin tz I_N\\
\cos z(1-t)\cos tz I_N& \cos z(1-t)\sin tz I_N\am .
\]
Let $c=\cos z, s=\sin z$ and $a=z(1-2t)$. Then substituting the
identity
$$
f=\ma (s+\sin a)I_N& (-c+\cos a)I_N\\
(c+\cos a)I_N& (s-\sin a)I_N\am=(s-cJ)+(\sin aJ_1-\cos a J_2)
=-Je^{zJ}+e^{aJ}J_2
$$
into \er{37} we obtain \er{36}.

Define the matrix $L={1\/2}(\wt M+\wt M^{-1})$ with
eigenvalues $\D_m={1\/2}(\t_m+\t_m^{-1}),m=1,2,...,N$, of
multiplicity two. The identity \er{TL-1} yields
\[
\lb{38}
L={\wt M+\wt M^{-1}\/2}=
{1\/2}\ma \vt(1,\cdot)+\vp'(1,\cdot)^\top & z(\vp(1,\cdot)-\vp(1,\cdot)^\top)\\
z^{-1}(\vt'(1,\cdot)-\vt'(1,\cdot)^\top) & \vt(1,\cdot)^\top+\vp'(1,\cdot)\am.
\]
Using \er{27} we obtain
\[
\lb{39} L(z)=L_1(z)+{L_2(z)\/2z^2}+{L_3(z)\/2z^3}+{O(e^{|\Im z|})\/z^{4}},\ \ \ \ L_1(z)=\cos z+{\sin z\/2z}V^0, 
\ \ \ \ |z|\to \iy,
\]
where
\[
\lb{310} L_2(z)=\int_0^1\!\! dt \int_0^t \sin z(t-s)
\ma a_{11}(t,s,z)& a_{12}(t,s,z)\\
a_{21}(t,s,z)& a_{22}(t,s,z)\am  ds, \ \ \ z\in \C,
\]
\[
\lb{311} a_{11}(t,s,z)=\sin z(1-t)\cos zs V(t)V(s)+ \cos
z(1-t)\sin zs V(s)V(t),\ \ \ a_{22}=a_{11}^\top,
\]
\[
\lb{312} a_{12}(t,s,z)=\sin z(1-t)\sin zs (V(t)V(s)-V(s)V(t)), \ \
\
\]
\[
\lb{313} a_{21}(t,s,z)=\cos z(1-t)\cos zs (V(t)V(s)-V(s)V(t)),
\]
and
\[
\lb{314} L_3(z)=\int_0^1\!\! dt\!\!  \int_0^t\!\! ds\!\!
\int_0^s\!\!   \sin z(t-s)\sin z(s-u)
\ma b_{11}(t,s,u,z)& b_{12}(t,s,u,z)\\
b_{21}(t,s,u,z)& b_{22}(t,s,u,z)\am  du,
\]
\[
\lb{315} b_{11}=\sin z(1-t)\cos zu V(t)V(s)V(u)+ \cos z(1-t)\sin
zu V(u)V(s)V(t),\ \ \ b_{22}=b_{11}^\top,
\]
\[
\lb{316} b_{12}=\sin z(1-t)\sin zu (V(t)V(s)V(u)-V(u)V(s)V(t)), \
\ \
\]
\[
\lb{317} b_{21}=\cos z(1-t)\cos zu (V(t)V(s)V(u)-V(u)V(s)V(t)).
\]
Recall that $ V^{0}=\!\!\int_0^1V(x)dx,\ \ B_{n}=\Tr
{(V^{0})^n\/n!}$.

\begin{lemma} \lb{T31}
For each $(r,V)\in \R_+\ts\mH$ asymptotics \er{1asD} and the following ones 
\[
\lb{318} 2^{2N}\det L(z)=\exp \rt(-2Niz+i{\Tr V^0\/z}
 +{i\|V\|^2+o(1)\/4z^3}\rt),
\]
\[
\lb{319} \Tr L_2(z)=\rt(-B_2+{i\|V\|^2+o(1)\/2z}\rt)\cos z, 
\]
\[
\lb{320} \Tr V_0L_2(z)=-3B_3\cos z+O(e^{|\Im z|}/z),
\]
\[
\lb{321} \Tr L_3(z)=-iB_3{\cos z\/2}+o(e^{|\Im z|}),
\]
hold as $y\ge r|x|, y\to \iy$. Moreover, the following identity  and  the  asymptotics are fulfilled
\[
\lb{322} 
L_2(\pi n)={(-1)^n\/4}\ma X^{(n)}+[\hat V^{cn},V^0] & [\hat V^{sn},V^0-\hat V^{cn}] \cr
[\hat V^{sn},V^0+\hat V^{cn}] & X^{(n)}-[\hat V^{cn},V^0]\am,
\]
\[
\lb{323} 
\D_m(z)=\cos \lt(z-{V_m^0\/2z}\rt)+{O(|\hat V^{(n)}|+n^{-1})\/n^2} \ \ \ as \ z=\pi n+O(1/n), \ m=1,..,N,
\]
where $X^{(n)}=-(V^0)^2+(\hat V^{cn})^2+(\hat V^{sn})^2$ and
 $[A,B]=AB-BA$ for matrix $A,B$.
\end{lemma}
\no {\it Proof.} 
Recall the simple fact: Let $A, B$ be matrices and 
 and $\s(B)$ be
spectra of $B$. If $A$ be normal, then $\dist\{\s(A),\s(A+B)\}\le |B|$
(see [Ka,p.291]).

 The diagonal operator $L_1(z)$ has the eigenvalues
$\D_m^0(z)=\cos z-V_{m}^0{\sin z \/2z}, m=1,..,N$ with the
multiplicity 2. Using the result from [Ka] and asymptotics
\er{39} we deduce that the eigenvalues $\D_m(z)$ of matrix $L(z)$
satisfy the asymptotics \er{1asD}.

 Using \er{310} we obtain
$$
\Tr L_2(z)=
2\int_0^1\int_0^t\!\!\!\sin z(t-s)\sin z(1-t+s)F_2(t,s)dtds
$$
$$
=\int_0^1\int_0^t\!\!\! (\cos z(1-t+s)-\cos
zm)F_2(t,s)dtds=I_1^0(z)- B_2\cos mz,
$$
since
$$
\int_0^1\int_0^t\!\!\! F_2(t,s)dtds={1\/2}\Tr
\lt(\int_0^1V(t)dt\rt)^2=B_2.
$$
Due to \er{61} we have $I_1^0(z)={i\|V\|^2+o(1)\/2z}\cos z$,
which yields \er{319}.

We show \er{320}. Let $G(t,s)=\Tr V_0(V(t)V(s)+V(s)V(t))$. Using
\er{310}, \er{61} we have
$$
\Tr V_0L_2(z)=\Tr V_0\int_0^1\int_0^t\!\!\! \sin z(t-s)z(1-t+s)
(V(t)V(s)+V(s)V(t))ds
$$
$$
={1\/2}\int^1dt\!\int_0^t\!\! (-\cos
z+\cos z(1-t+s))G(t,s)ds=-3B_3\cos z +O(e^{|\Im z|}/z),
$$
since
$$
\int_0^1dt\int_0^t\!\!\! \Tr V_0(V(t)V(s)+V(s)V(t))ds
=\Tr V_0\lt(\int_0^1V(t)dt\rt)^2=6B_3.
$$
Consider $\Tr L_3$. The identity \er{314} gives
$$
\Tr L_3(z)=\int_0^1dt\int_0^t\!\! ds\! \int_0^s\!\!\!
\sin z(t-s)\sin z(s-u)
\Tr (b_{11}+b_{22})du
$$
$$
=\int_0^1dt\int_0^t\!\!ds\! \int_0^s\!\!
\sin z(t-s)\sin z(s-u)\sin z(1-t+u)Rdu,
$$
where $R=\Tr (V(t)V(s)V(u)+V(u)V(s)V(t))$. Using the identity
$$
4\sin z(t-s)\sin z(s-u)\sin z(1-t+u)=-\sin z+P,
$$$$P=\sin z(1-2s+2u)-\sin z(1-2t+2s)-\sin z(1-2t+2u),
$$
we get
$$
L_3=-L_3^0{\sin z\/2} +L_3^1,\ \ \
L_3^0={1\/2}\int_0^1dt\int_0^t\! ds\!\int_0^s\! Rdu,\ \
L_3^1={1\/4}\int_0^1dt\int_0^t\! ds\!\int_0^s\! P Rdu,
$$
where
$$
L_3^0={\Tr\/2}\int_0^1V(t)dt\int_0^t\! ds\!\int_0^s\!
(V(s)V(u)+V(u)V(s))du=
{\Tr\/2}\int_0^1V(t)\rt(\int_0^t\! V(s) ds\rt)^2
=B_3.
$$
Due to \er{63} we obtain $L_3^1=o(e^{|\Im z|})$, which yields \er{321}.

We will show \er{318}. Asymptotics \er{39} yields
\[
\lb{324}
{L\/\cos z}=I_{2N}+S,\ \
S=i{V^0I_{2N}\/2z}+{L_2\/2z^2\cos z}+{L_3\/2z^3\cos z}+O(z^{-4}),\ \
\]
as $|z|\to \iy, y\ge r|x|$,
since $\tan z=i+O(e^{-y})$. In order to use the identity
$$
\det (I+S)=\exp \rt(\Tr S-\Tr {S^2\/2}+\Tr {S^3\/3}+..\rt),\ \
|S|\to 0,
$$
we need the traces of $S^m, m=1,2,3$. Due to \er{319}-\er{322} we
get
\[
\lb{325}
{\Tr S^3\/3}=-i\Tr {(V^0)^3I_{2N}\/3(2z)^3}+O(z^{-4})
=-i {B_3\/2z^3}+O(z^{-4}),
\]
\[
\lb{326}
-\Tr{S^2\/2}={\Tr\/2}\rt({(V^0)^2\/(2z)^2}I_{2N}-
2i{V^0L_2\/4z^3\cos z}+O(z^{-4})\rt)=
{B_2\/2z^2}+i{3B_3\/4z^3}+O(z^{-4})
\]
and
\[
\lb{327}
\Tr S=\Tr \rt(i{V^0\/2z}+{L_2\/2z^2\cos z}+
{L_3+O(z^{-1})\/2z^3\cos z}\rt)
=i{B_1\/z}+\rt(-{B_2\/2z^2}+i{\|V\|^2\/4z^3}\rt)-i{B_3+o(1)\/4z^3}
\]
and summing \er{325}-\er{327} we get \er{318}. 

We will show \er{322}. Let $z=\pi n, c_t=\cos \pi nt, s_t=\sin \pi nt,$ and $K_{ts}^\pm=V(t)V(s)\pm V(s)V(t)$. Using \er{311}, we have
$$
4(-1)^n\sin \pi n(t-s)a_{11}(t,s,\pi n)=4(s_tc_s-s_sc_t)(-s_tc_sV(t)V(s)+s_sc_tV(s)V(t))
$$
$$
=((c_{2t}-1)(1+c_{2s})+s_{2s}s_{2t})V(t)V(s)+
((c_{2t}+1)(1-c_{2s})+s_{2s}s_{2t})V(s)V(t)
$$
$$
=(-1+c_{2t}c_{2s}+s_{2t}s_{2s})K_{ts}^++((c_{2t}-c_{2s})K_{ts}^-
$$
and the integration yields
\[\lb{328}
\int_0^1\!\!\!\!  dt\!\!  \int_0^t \!\! \sin z(t-s)
a_{11}(t,s,z)ds={1\/2}\int_0^1\!\!  \!\! \int_0^1 \sin z(t-s)
a_{11}(t,s,z)dtds
=
{(-1)^n\/4}\rt(X^{(n)}+[\hat V^{cn},V^0]\rt),
\]
where $X^{(n)}=-(V^0)^2+(\hat V^{cn})^2+(\hat V^{sn})^2$. Using \er{312}, we obtain 
$$
(-1)^n\sin \pi n(t-s)a_{12}(t,s,\pi n)=(s_tc_s-s_sc_t)(-s_ts_s)K_{ts}^-=(s_{2t}-s_{2s}+s_{2s}c_{2t}-
c_{2s}s_{2t}){K_{ts}^-\/4}
$$
and the integration yields
\[
\lb{329}
\int_0^1\!\! dt \int_0^t \sin z(t-s)
a_{12}(t,s,z)ds={(-1)^{n}\/4}\rt(\hat V^{sn}V^0-V^0\hat V^{sn}-
\hat V^{sn}\hat V^{cn}+\hat V^{cn}\hat V^{sn}\rt).
\]
The similar arguments give
\[
\lb{330}
\int_0^1\!\! dt \int_0^t \sin z(t-s)
a_{21}(t,s,z)ds={(-1)^{n}\/4}\rt(V^0\hat V^{sn}-\hat V^{sn}V^0-
\hat V^{sn}\hat V^{cn}+\hat V^{cn}\hat V^{sn}\rt).
\]
Combining the identities \er{328}-\er{330} and using $a_{22}=a_{11}^\top$ we obtain \er{322}.

The diagonal operator $L_1(z)$ has the eigenvalues
$\D_m^0(z)=\cos z-V_{m}^0{\sin z \/2z}, m=1,..,N$ with the
multiplicity 2. Using the result from [Ka] and asymptotics
\er{322} we deduce that the eigenvalues $\D_m(z)$ of matrix $L(z)$
satisfy the asymptotics \er{323} as $z=\pi n+O(1/n)$.
$\BBox$

Recall that $D^0(\pm 1,z)= (e^{iz}\mp 1)^N(e^{-iz}\mp 1)^N=
2^N(1\mp\cos z)^N$.

\begin{lemma} \lb{T32}
Let $V\in \mH$ and let $A= e^{|\Im z|+\vk},\ 
\vk={\|V\|\/|z|_1},\  |z|_1=\max\{1,|z|\}$. Then

\no i) The following estimates are fulfilled:
\[
\lb{331} |\x_m(z)|\le (2NA)^m,\ \ \ |\x_m(z)-\x_m^0(z)|\le 2\vk
(2NA)^m,\ \ \ m=1,..,N.
\]
\[
\lb{332}
 |D(\pm 1,z)-D^0(\pm 1,z)|<C_N\vk A^N,\ \ \ \ C_N=4(2N)^N.
\]
\no ii) For each integer $n_0>4^NC_N\|V\|$ the function $D(1,z)$ has
exactly $N(2n_0+1)$ roots, counted with multiplicity, in the disc
$\{|z|<\pi(2n_0+1)\}$ and for each $|n|>n_0$, exactly $2N$ roots,
counted with multiplicity, in the domain $\{|z-2\pi
n|<{\pi\/2}\}$. There are no other roots.

\no iii) For each integer $n_0>4^NC_N\|V\|$ the function $D(-1,\l)$
has exactly $2Nn_0$ roots, counted with multiplicity, in the disc
$\{|z|<2\pi n_0\}$ and for each $|n|>n_0$, exactly $2N$ roots, counted
with multiplicity, in the domain $\{|z-\pi (2n+1)|<{\pi\/2}\}$.
There are no other roots.

\no iv) Assume that $V_i^0\ne V_j^0$ for all $i\ne j\in \o_s$ for some $s=1,..,N_0$. Then there exists 
integer $n_0\ge 1$ such that the function $\r_s$
has exactly $2N_s(N_s-1)n_0$ roots, counted with multiplicity, in the disc $\{|z|<\pi (n_0+{1\/2})\}$ and for each $|n|>n_0$, exactly $N_s(N_s-1)$ roots, counted
with multiplicity, in the domain $\{|z-\pi n|<{\pi\/2}\}$.
There are no other roots.
\end{lemma}

\no {\it Proof.} We prove \er{331} by induction. Let the first estimate in \er{331} hold for $\x_m$. Using \er{29} we obtain $|\x_m(z)|\le \m^m, \m=2NA$ for $m=0,1,2$. If $m\ge 2$, then substituting the estimate
$|\x_j(z)|\le \m^j$ and $|T_j(z)|\le A^j$ (see \er{29})  into \er{15} we have
$$
|\x_{m+1}(z)|\le{2N\/m+1}\sum_0^{m}A^{m+1-j} \m^{j}\le
{2NA^{m+1}( (2N)^{m+1}-1)\/(m+1)(2N-1)}\le \m^{m+1}.
$$
We shall show the second estimate in \er{331}. Let
$p_n=|\x_n(z)-\x_n^0(z)|$. The recurrent identities \er{15} and
$|T_j(z)-T_j^{(0)}(z)|\le \vk jA^j, |\x_j(z)|\le (2NA)^j $ and \er{29} give
$$
p_n\le{2N\/n}\sum_0^{n-1}|T_{n-k}(z)\x_k(z)-T_{n-k}^{(0)}(z)\x_k^0(z)|
 \le
{2N\/n}\sum_0^{n-1}\lt(|\x_k(z)||T_{n-k}(z)-T_{n-k}^{(0)}(z)|
+|T_{n-k}(z)|p_k\rt)
$$
$$
\le {2N\/n}\sum_1^{n-1}\m^k\vk (n-k)A^{n-k}+A^{n-k}\vk k  \m^k
=2NA^n\vk\sum_1^{n-1}(2N)^k\le 2\vk \m^n
$$
where we used $p_1\le \vk A$, see \er{29}. Thus we get \er{331}.

 ii) We have $ D(1,\cdot)=\x_N+2\sum_0^{N-1}\x_n, $
where $\x_n$ are given by \er{15}. Recall that $\x_n^0=\x_n$ at $V=0$. Hence using \er{331} and $\m=2NA$ we have
\[
\lb{333}
 |D(1,z)-D^0(1,z)|\le |\x_N(z)-\x_N^0(z)|+
 2 \sum_1^{N-1}|\x_n(z)-\x_n^0(z)|
 \le 2\vk \m^N+2\sum_1^{N-1}2\vk \m^n.
\]
Thus  $|D(1,z)-D^0(1,z)|<2\vk \m^N(1+{2\/\m-1})<4\vk \m^N
$ since $\m\ge 4$ which yields \er{332}.

  Let $n_1>n_0$ be another integer.
Introduce the contour $C_n(r)=\{z:|z-\pi n|=\pi r\}$. Consider the
contours $C_0(2n_0+1),C_0(2n_1+1),C_{2n}({1\/4}),|n|>n_0$. Note that
$\vk\le {1\/4^4}$ on all contours. Then \er{327} and the estimate
$e^{{1\/2}|\Im z|}<4|\sin{z\/2}|$ on all contours yield 
\[
\lb{334}
\lt|D(1,\l)-\lt(2\sin{z\/2}\rt)^{2N}\rt|\le C_N\vk A^N=
 C e^{N|\Im z|}\le
(C4^N) \lt|2\sin{z\/2}\rt|^{2N}<{1\/2}\lt|2\sin{z\/2}\rt|^{2N}
\]
where $C=|z|^{-1}C_N\|V\|e^{N\vk}$.
Hence, by Rouch\'e's theorem, $D(1,z)$ has as many roots, counted
with multiplicities, as $\sin^{2N}{z\/2}$ in each of the bounded
domains and the remaining unbounded domain. Since
$\sin^{2N}{z\/2}$ has exactly one root of the multiplicity $2N$ at
$2\pi n$, and since $n_1>n_0$ can be chosen arbitrarily large, the
point ii) follows.

 iii) The proof for $D(-1,\l)$ is similar.
 
iv) Consider the case $N_0=1$, the proof for $N_0\ge 2$ is similar.  Asymptotics \er{1asD} yields
$\D_j(z)-\D_m(z)=(V_{j}^0-V_{m}^0){\sin z
\/2z}+O(z^{-2}e^{|\Im z|}),\ \ |z|\to\iy.
$
Then
\[
\lb{335}
\r(z)=\prod_{j<m}
(\D_j(z)-\D_m(z))^2=c_0\rt({\sin z\/2z}+O({e^{|\Im z|}\/z^{2}})\rt)^{N(N-1)}\!\!\!\!,\ \ 
c_0=\!\!\prod _{j<m} \!\!(V_{j}^0-V_{m}^0)^2.
\]
 Let $n_1>n_0$ be another integer.
Introduce the contour $C_n(r)=\{z:|z-\pi n|=\pi r\}$. Consider the
contours $C_0(2n_0+1),C_0(2n_1+1),C_{2n}({1\/4}),|n|>n_0$. Note that
$\vk\le {1\/4^4}$ on all contours. Then \er{335} and the estimate
$e^{|\Im z|}<4|\sin z|$ on all contours (for large $n_0$) yield
$$
\r(z)=\r^0(z)(1+O(z^{-1})),\ \ \ \r^0(z)=c_0\rt({\sin z
\/2z}\rt)^{N(N-1)}.
$$
Hence, by Rouch\'e's theorem, $\r$ has as many roots, counted
with multiplicities, as $\r^0$ in each of the bounded domains and
the remaining unbounded domain. Since $\r^0$ has exactly one root
of the multiplicity $N(N-1)$ at $\pi n\neq 0$, and since $n_1>n_0$
can be chosen arbitrarily large, the point iv) follows.
$\BBox$

{\bf Proof of Theorem   \ref{T2}} i)
 We determine asymptotics \er{T2-1} for $z_{m,m}^{n\pm}=\sqrt{\l_{m,m}^{n\pm}}$ as $n\to \iy$. Lemma \ref{T32} yields
$|z_{m,m}^{n\pm}-\pi n|<{\pi\/2}$ as $n\to \iy,m=1,2,..,N$. 
Lemma \ref{T31} gives 
$
\D_m(z)=\cos (z-{V_m^0\/2z})+O(1/z^2)$ as   $z=\pi n+O(1)$.
Using the identity 
$\D_j(z_{m,m}^{n\pm})=(-1)^n$, we have $z_{m,m}^{n\pm}=\pi n+O(1/n)$. We shall determine sharper asymptotics.

Recall that the modified monodromy matrix $\wt M(z)$ is given by \er{35} and it has
the same eigenvalues as $M(z)$. 
Define the local parameter $\m$ by $z=\pi n+\ve\m, \ve={1\/2\pi n}$.
Note that $\l=z^2=(\pi n)^2+\m+(\ve\m)^2$.
Asymptotics  \er{36} gives 
\[
\lb{336} (-1)^n\wt M(z)=e^{\ve \m J}\rt(I_{2N}-\ve J(V_0+
J\hat V(\pi n)J_2)+O(\ve)\rt)=I_{2N}-\ve J(V_0+
\hat \cV^n+O(\ve)-\m),\
\ \ \ 
\]
where
$$
\hat \cV^{n}=J\hat V(\pi  n)J_2=J\int_0^1\!\! V(t)e^{-2tzJ}dtJ_2=J(\hat V^{nc}-J\hat V^{ns})J_2=J_1\hat V^{nc}+J_2\hat V^{ns}.
$$
Hence we have the asymptotics \er{T2-1}, since
$$
0=\det \lt((-1)^n\wt M(z)-I_{2N}\rt)=\det (-\ve J)\det \rt(V_0+
\hat \cV^n+O(\ve)-\m\rt).
$$

Consider the case $V_{1}^0<...<V_{N}^0$. We shall determine  asymptotics \er{T2-4} for the case $\a=(m,m), m=1,..,N$.
Let $z_{m,m}^{n\pm}=\pi n+\ve\m$ and $\m-V_{1}^{0}=\x\to 0$.
Using the simple transfomation (unitary), i.e., changing
the lines and columns, we obtain
$$
\det (V^0+\hat \cV^{(n)}+O(\ve)-\m)=
\det \ma A_1-\x& A_2\\ A_{3}& A_4-\x \am=\det (A_4-\x)\det K(\x),
$$$$
 A_1=
\ma \hat V_{11}^{cn}& \hat V_{11}^{sn}\\ \hat V_{11}^{sn}& -\hat V_{11}^{cn}\am
+O(\ve),\ \ \ \ 
A_4=\diag \{V_{2}^0-V_{1}^0,...,V_{N}^0-V_{1}^{0}\}+A_5,\ \ 
$$
$$
A_2,A_3, A_5=O(\d_n),\ \ \d_n=|\hat V^n|+\ve,\  \ \
K(\x)=A_1-\x-A_2(A_4-\x)^{-1}A_3+O(\ve)
$$
 We have $K(\x)=A_1-\x+O(\f), \f=\ve+|\hat V^n|^2 $, which yields
$$
0=\det K(\x)=\x^2-|\hat V_{11}^n|^2-\x b_1+\hat V_{11}^{cn}b_2+\hat V_{11}^{sn}b_3+O(\f^2), \ \ \ b_1,b_2,,b_3=O(\f)
$$
where $b_1,b_2,b_3$ are analytic functions of $\x$. Rewriting
the last equation in the form $(\x+\a)^2=(|\hat V_{11}^n|+\b)^2+O(\f^2), \a,\b=O(\f)$ and using the estimate $\sqrt{x^2+y^2}-x\le y$ for $x,y\ge 0$ we get 
$\x=\pm |\hat V_{11}^n|+O(\f)$, which yields \er{T2-4}
for the case $\a=(m,m)$.

Consider the resonances.
Assume that $V_i^0\ne V_j^0$ for all $i\ne j\in \o_s$ for some 
$s=1,..,N_0$. By Lemma \ref{T32}, the zeros of $\r_s$ have the form
$z_\a^{n\pm}, \a=(j,j'), j,j'\in \o_s, j<j', n\in \Z\sm \{0\}$
and satisfy $|z_\a^{n\pm}-\pi n|<\pi/2$. 

Asymptotics \er{1asD} yields
$\D_j(z)-\D_{j'}(z)=(V_{j}^0-V_{j'}^0){\sin z
\/2z}+O(z^{-2}e^{|\Im z|}),\ \ |z|\to\iy$. Then 
$   |z_\a^{n\pm}-\pi n|<\pi/2$ yields $z_\a^{n\pm}=\pi n+O(1/n)$ as $n\to \iy$.

We have
the identity   $\D_j(z)-\D_{j'}(z)=0$ at $z=z_\a^{n\pm}$.
Then using \er{323} we have
$$
\cos\lt(z_\a^{n\pm}-{V_j^0\/2\pi n}\rt)- \cos\lt(z_\a^{n\pm}-{V_{j'}^0\/2\pi n}\rt)=2(-1)^n\sin {V_{j'}^0-V_j^0\/4\pi n}\sin\lt(z_\a^{n\pm}-\pi
n-{V_j^0+V_{j'}^0\/4\pi n}\rt)={O(\d_n)\/n^2}
$$
which yields \er{T2-2}, i.e.,
\[
z_\a^{n\pm}=\pi n+\ve(a_++\wt z_\a^{n\pm}),\ \
a_\pm={V_j^0\pm V_{j'}^0\/2},\ \ \ \ve={1\/2\pi n},\ \  \   \\wt z_\a^{n\pm}=\d_n.
\]
We shall show that for large $n$ in the neighborhood of each
$\pi n+\ve a_+$ the function $(\D_j(z)-\D_{j'}(z))^2$
has two real zeros resonances (counted with multiplicity).
Introduce the functions
\[
\lb{338}
f_m(\m)=2(2\pi n)^2(1-(-1)^n\D_m(\pi n+\ve \m))=(\m-V_m^0)^2+
O(\d_n),\ \ \ \d_n=|\hat V^{(n)}|+\ve.
\]
For the case $\m\to a_+$ we get
\[
\lb{339}
f_m(\m)=(a_+-V_m^0)^2+o(1),  \ \ m=1,..,N, \ \ and \  \  
f_m(\m)=a_-^2+o(1),\ \ m=j,j'. 
\]
Hence the function $f_j-f_{j'}$ (maybe)
has the zeros, but the functions $f_j-f_m, m\ne j,j'$
have not zeros in the neighborhood of the point $a_+$.

Note that these functions are real outside the small 
neighborhood of $a_+$, otherwise for any complex branches
there exists a complex conjugate branch, but the asymptotics
\er{338} show that such branches are absent.

We have two cases: (1) let $f_m(\m),m=j,j'$
be real in some small neighborhood of $a_+$. 
Then the function $f_j-f_{j'}$ has at least
one real zero, since by Theorem \ref{T1}, the functions $f_j, f_{j'}$ are strongly monotone. Thus $(f_j-f_{j'})^2$ has at least 2 real zeros.

(2) Let $f_m(\m),m=j,j'$
be complex in some small neighborhood of $a_+$. Then they have
at least two real branch  points. Thus $(f_j-f_{j'})^2$ has at least 2 real zeros.

Hence $(f_j-f_{j'})^2$ has exactly two real zeros, since
the number of resonances (in the neighborhood of the point
$\pi n$) is equal to $N_s(N_s-1)$.

We determine the sharp asymptotics of resonancees.
Define the unitary matrix $P={1\/\sqrt 2}(J_1+iJ)=P^*, P^2=I_{2N}$.
Using the identities 
$$
PJP=-iJ_1,\ \ \ PJ_1P=iJ,\ \ \ PJ_2P=-J_2,\ \ \ 
$$
$$
\mV_{n}=-iJ_1P\hat \cV^n P=-iJ_1(iJ\hat V^{cn}-J_2\hat V^{sn})=J_2\hat V^{cn}+iJ\hat V^{sn}=\ma 0&\hat V^{(n)} \\ ({\hat V^{(n)}})^* &0\\ \am
$$
we have
$$
A=P\rt( {(-1)^n\wt M(z)-I\/i\ve}\rt)P=J_1(V^0+P\hat V^n P+O(\ve)-\m)=
J_1(V^0-\m)+i\mV_{n}+O(\ve).
$$
The operator $A-a_-$ has the eigenvalue $\x_0={(-1)^n\t_{n,s}-1\/i\ve}-a_-$
of multiplicity two, since $\t_{n,s}=\t(z_\a^{n\pm})=(-1)^ne^{i\ve(a_-+o(1))}$.
The operator $J_1(V^0-a_+)-a_-=\ma V^0-V_j^0& 0\\ 0& -V^0+V_{j'}^0\am$ has two eigenvalue ($= 0$) and other eigenvalues
are not zeros.
Using the simple transformation (unitary), i.e., changing
the lines and columns, $\m=a_++r\in \R, r=\wt z_\a^{n\pm}\to 0, v=\hat V_{jj'}^{(n)},$ we obtain
$$
F(\x)=\det (A-a_--\x )=
\det \ma A_1-\x& A_2\\ A_{4}& A_3-\x \am=\det (A_3-\x)\det \rt(K(\x)-\x I_2\rt),\ \
$$$$
 A_1=
\ma -r& iv\\ i\ol v& r\am +O(\ve),\ \ \  \  K(\x)=A_1-A_2(A_3-\x)^{-1}A_4=\ma -r+a_1& iv+a_4\\ i\ol v+a_3& r+a_2\am ,
$$$$
A_3=\diag \{V_m^0-V_j^0-r, m\neq j\}\os
\diag \{V_{j'}^0-V_m^0+r, m\neq j'\}+A_4,\ \ A_2,A-3, A_4=O(\d_n),
$$
the function $a_1,a_2,a_3,a_4=O(\f), \f=\ve+|\hat V^{(n)}|^2$ and they analytic with respect to $\x$ in some small disk.
The function $F=\det \rt(K(\x)-\x I_2\rt)$ has the form 
\[
F(\x)=\x^2-r^2+|v|^2+a_1(r-\x)+a_2(-r-\x)-iva_3-ia_4\ol v-a_4a_3
=(\x-\x_0)^2(1+O(\x-\x_0))
\]
 for $\x\to 0$ where $\x_0={((-1)^n\t_{n,s}-1)\/i\ve}-a_-$ is the zero of $F$ of multiplicity two.
Then $\x_0=O(\f)$ and we have $(r-\a)^2=(|v|-\b)^2+O(\f^2)$ where $\a,\b=O(\f)$. Then using the estimate $\sqrt{x^2+y^2}-x\le y$ for $x,y\ge 0$ we get $r=\pm |v|+O(\f)$.
$\BBox$

{\bf Proof of Corollary \ref{Tc1}} (i) Let $N_G=\iy$. Then, due to the Lyapunov Theorem, Theorem \ref{T2},
there exists a real sequence $\l_k\to+\iy$ as $k\to\iy$, such that
$\l_k\in\g_{n_k}^{\{j(m),m\}}$ for each $m=1,..,N$. Hence,
$\cap_{m=1}^N\g_{n_k}^{\{j(m),m\}}\ne \es$. 
Using  asymptotics \er{T2-4} and $k\to \iy$, we obtain $V_1^0+V_{j(1)}^0=...=V_N^0+V_{j(N)}^0$. Moreover, the estimates $V_1^0<...<V_N^0$ yield
$V_{j(1)}^0>...>V_{j(N)}^0$, i.e. $j(1)=N$, $j(2)=N-1$, ... Then,
$V_1^0+V_N^0=V_2^0+V_{N-1}^0=...$, which give a contradiction.\\
(ii) Let $2a=V_1^0+V_N^0=V_2^0+V_{N-1}^0=..$. Due to \er{T2-4}, $(\pi n_k)^2+a\in \cap_{m=1}^N\g_{n_k}^{\{N+1-m,m\}}$ as $k\to\infty$. Then the Lyapunov
Theorem yields $(\pi n_k)^2+a\notin\s(\mL)$, $k\to\iy$, i.e. $N_G=\iy$.
$\BBox$

\section {Harmonic functions}
\setcounter{equation}{0}

In this Sect. we will prove Theorems \ref{T41} and \ref{T42}
about the properties of the quasimomentum.
Recall that the Lyapunov function $\wt\D_s(\z)$ is analytic on some 
 $N_s$--sheeted Riemann surface $\mR_s$ and $\mR=\cup_1^{N_0} \mR_s$.
Let $z=x+iy\in\C$ be the natural projection of $\z\in\mR$, $\cB_\D$ be the set of all branch points of the Lyapunov function and $\mR^\pm=\{\z\in\mR:\pm\Im\z>0\}$. We define the simply
connected domains $\mR_0^\pm\ss\C_\pm$ and a domain $\mR_0$  by
$$
\mR_0^\pm=\C_\pm\sm\b_\pm,\ \ \ \b_\pm={\bigcup}_{\b\in\cB_\D\cap\C_\pm}[\b,\b\pm i\infty).
$$
$$
\mR_0=\C\sm \rt(\b_+\cup \b_-\cup \b_0\rt),\ \
\b_0=\{z\in \R: \D_m(z)\notin \R \ for \ some \ m\in\{1,..,N\}\}
$$
Due to the Lyapunov Theorem, $\D(\z)\notin [-1,1]$, $\z\in\mR^+$. Recall that
$q(\z)=|\log\e(\D(\z))|$ is the single-valued on $\mR^+$ imaginary part of the (in general, many-valued on $\mR^+$) quasimomentum $k(\z)=p(\z)+iq(\z)=\arccos\D(\z)=i\log\e(\D(\z))$,
where
$$
\e(z)\equiv z+\sqrt{z^2-1},\ \ \ \ \ \e:\C\sm [-1,1]\to \{z\in \C: |z|>1\}.
$$
We denote by $q_m(z)$, $z\in\C_+$, $m=1,..,N$, the branches of $q(\z)$ and by $p_m(z)$,
$k_m(z)$, $z\in\mR_0^+$, the single-valued branches of $p(\z)$, $k(\z)$, respectively.

\begin{theorem} \lb{T41}
Assume that  $V\in \mH$  and fix some $j=1,..,N_0$ . Then
the function $\wt q_s(\z)=\log |\e(\wt\D_s(\z))|$
 is subharmonic on the Riemann surface $\mR_s$ 
and the following asymptotics are fulfilled:
\[
\lb{41} \wt q_s(\z)=y+O(1/|z|),\ \ \ \ y>r|x|,\ \ \ {\rm 
\ any}\  \ r>0, \ \
\]
\[
\lb{42} \wt q_s(\z)=y+O(|z|^{-{1\/2}}),\ \ \ \ 
\]
as $|\z|\to \iy, \  \z\in \mR_s$.
Moreover, let $\D_j$ be analytic on some bounded interval
$Y=(\a,\b)\ss \R$ for some $j\in \o_s$. Then

\no i) If $\D_j(z)\in \R\sm [-1,1]$ for all $z\in Y$, then $k_j(\cdot)$ has an analytic
extension from $\mR_0^+$ into $\mR_0^+\cup\mR_0^-\cup Y$ such
that
\[
\lb{43} \Re k_j(z)={\rm const}\in \pi \Z,\ \
z\in Y,
\]
\[
\lb{44} q_j(z)=q_j(\ol z)>0,\ \ \ \ z\in \mR_0^+\cup
\mR_0^-\cup Y.
\]
\no ii) If we assume $\D_j(z)\notin \R$ for 
any $z\in Y$, then there exists a
branch $\D_i,i\in \o_s$ such that $\ol\D_i(z)=\D_j(z)$ for any
$z\in Y$. The functions $\D_j(z)$ and $k_i+k_j$ have
analytic extensions from $\mR_0^+$ into
$\mR_0^+\cup\mR_0^-\cup Y$ such that
\[
\lb{45}
\D_j(z)=
\cases {\D_j(z)\ \ \ &if\ \ \ \ $z\in \mR_0^+$\cr
       \ol\D_i(\ol z)\ \ \ &if\ \ \ \ $z\in \mR_0^-$\cr},     
\]
\[
\lb{46}
 p_j(z)+p_i(z)={\rm const}\in 2\pi \Z,\ \
z\in Y,\ \ \
q_j(z)=q_i(z),\ \ \ \ \ z\in Y,
\]
\[
\lb{47}
 q_j(z)+q_i(z)=q_j(\ol z)+q_i(\ol z)>0,\ \
\ z\in \mR_0^+\cup \mR_0^-\cup Y.
\]
\end{theorem}
\no {\it Proof.}  By Theorem \ref{T1}, the function $\wt\D_s$ is analytic on $\mR_s$ and the function $\e(\cdot)$ is subharmonic on $\C$. Then $\wt q_s(\z)=\log |\e(\wt\D_s(\z))|$ is subharmonic on $\mR_s$.  Using the asymptotics \er{dx}, \er{1asD} we obtain
$$
\wt q_s(\z)=\log |2\cos z+O(y^{-1}e^y)|=y+O(y^{-1}),\ \ \ |z|\to
\iy, \ y>r|x|,
$$
which yields \er{41}. Due to \er{1asD} and Lemma \ref{TA3}, we have
$$
\wt q_s(\z)=\log |\e(\cos z+O(e^{|\Im z|}/z))|=
\log |\e(\cos z)|+O(1/\sqrt z)=y+O(1/\sqrt z))
$$
which yields \er{42}.

i)  Due to $\D_j(z)=\cos k_j(z)$ we obtain \er{43}. The
real part of $k_j$ is a constant on $Y$, then $k_j$ has
an analytic extension from $\mR_0^+$ into $\mR_0^+\cup\mR_0^-\cup
Y$. Moreover, $q_j$ has an harmonic extension from $\mR_0^+$
into $\mR_0^+\cup\mR_0^-\cup Y$ by $q_j(z)=q_j(\ol z),\ \
z\in \mR_0^\pm$.

ii) By Theorem \ref{T1}, each polynomial $\F_s(\n,z)=\prod_{n\in \o_s} (\n -\D_n(z)),
z^2\in \R$ is real for $\n\in \R$. Then for $\D_j$ there exists a
$\D_i$ such that $ \D_j(x)=\ol \D_i(x), x\in Y$.
Then by the Morer Theorem, the function $\D_j$ has an analytic
extension  given by \er{45} from $\mR_0^+$ into $\mR_0^+\cup
\mR_0^-\cup Y$. Using $\D_m(z)=\cos k_m(z), m=j,i$  we obtain $0=\D_j(x)-\ol\D_i(x)=-2\sin {k_j(x)+\ol
k_i(x)\/2} \sin {k_j(x)-\ol k_i(x)\/2}$. Thus we get \er{46},
\er{47} since $q_j(x)>0, q_i(x)>0$ on $Y$.
$\BBox$

Recall the needed properties of the functions $v\in\cS\cC$
defined in Sect. 1 and $w=u+iv$. It is well known, that $u\in C(\ol\C_+)$ and
 ${1\/2\pi}\D v=\m_v$ (in a sense of
distribution) is a so-called Riesz measure of the function $v$.
Moreover, the following identities are fulfilled:
\[
\lb{48} \pi\m_v((x_1,x_2))=u(x_2)-u(x_1), \ \ \ \ {\rm for \ any
}\ \ x_1<x_2,\  x_1,x_2\in \R,
\]
\[
\lb{49} {\pa v(z)\/\pa y}=y\int_{\R}{d\m_v(t)\/(t-x)^2+y^2},\ \ \
z=x+iy\in\C_+,
\]
which yields ${\pa v(z)\/\pa y}\ge 0, z\in\C_+$. Moreover,
$v(x)=v(x\pm i0), x\in \R.$ It is well known that if $v\in
\cS\cC$, then
\[
\lb{410} \int_{\R}{d\m_v(t)\/1+t^2}<+\iy,\ \ \
\lim\sup_{z\to\iy}{v(z)\/|z|}=
\lim_{y\to+\iy}{w(iy)\/iy}=\lim_{y\to+\iy}{v(iy)\/y}=\lim_{x\to
\pm \iy}{u(x)\/x}\ge 0.
\]
 Now we recall the well known fact (see [Ah]).

\no {\bf Theorem (Nevanlinna).}
 {\it \no i) Let $\m$ be a Borel measure on
$\R$ such that  $\int_{\R}(1+x^{2p})d\m(x)<+\iy$ for some
$p\in\Z_+$. Then for each $r>0$ the following asymptotics is valid
$$
\int_{\R}{d\m(t)\/t-z}=-\sum_{k=0}^{2p}{Q_k\/z^{k+1}}+
{o(1)\/z^{2p+1}},\ \ \ |z|\to\iy,\ y>r|x|,\ \ \ \
Q_n=\int_{\R}x^nd\m(x),\ 0\le n\le 2p.
$$
\no ii) Let $F$ be an analytic function in $\C_+$ such that $\Im
F(z)\ge 0,\ z\in\C_+$ and
\[
\lb{411}
\Im F(iy)=c_0y^{-1}\!+\!\dots\!+\!c_{2p-1}y^{-2p}\!+\!O(y^{-2p-1})
\ \ \ {\rm as}\ \ y\!\to\!\iy
\]
for some $c_0,..c_{2p-1}$ and $p\ge 0$.  Then $
F(z)=C+\int_{\R}{d\m(t)\/t-z},\ z\in\C_+$, for some Borel measure
$\m$ on $\R$ such that $\int_{\R}(1+x^{2p})d\m(x)<\iy$ and $C\in
\R$.}

\begin{theorem} \lb{T42}
Assume that $V\in \mH$. Then
the function $w={1\/N}\sum_1^N k_m$ is analytic in $\C_+$;
the function $v=\Im w={1\/N}\sum_1^N q_m$ belongs to $\cS\cK_0^+\cap
C(\C)$ and it is positive harmonic in $\O=\C_+\cup \C_-\cup g$.
Moreover, the following asymptotics and identities hold
\[
\lb{412} u(z)={\rm const}\in {\pi\/N} \Z,\ \ \ \ \ v(z)>0,\ \ \ \ \
z\in g_n,
\]
\[
\lb{413} v(z)=y+O(1/y),\ \ \ \ |z|\to \iy,\ \ y>r|x|,
\ \ \ {\rm \ any}\  \ r>0,
\]
\[
\lb{414} v(z)=y+O(|z|^{-{1\/2}})\ \ \ \ \ |z|\to\iy.
\]
\end{theorem}
\no {\it Proof.}  Recall that $k_m(z)=i\log\e(\D_m(z))$, $z\in\mR_0^+$, where $\D_m(z)$, $m=1,..,N$, are the
branches of the function $\D$ analytic on $\mR$ . Since $\D(\z)\notin [-1,1]$,
$\z\in\mR^+$, the function $\e(\D(\z))$ is analytic on $\mR^+$ and $|\e(\D(\z))|>1$,
$\z\in\mR^+$. Recall that $\cB_\D$ is the set of all branch points of the function $\D$. Define the function
$
F(z)\equiv\prod_{m=1}^N\e(\D_m(z)),\  z\in\C_+.
$
Due to the symmetry of $F$ with respect to the permutations of $\D_1,..,\D_n$, we deduce that
$F$ is analytic in the domain $B(z_0,r)\sm\{z_0\}$ for any $z_0\in\C_+$ and some $r=r(z_0)$,
including the case $z_0\in \C_+\cap\cB_\D$. Note that $F$ is bounded in $B(z_0,r)$. This
implies  $F$ is analytic in $B(z_0,r)$ for any $z_0\in\C_+$, i.e. $F$ is analytic in $\C_+$
and $|F(z)|>1$, $z\in\C_+$. Therefore, the averaged quasimomentum
\[\lb{417}
w(z)={i\/N}\log F(z)={1\/N}\sum_{m=1}^N i\log\e(\D_m(z))={1\/N}\sum_{m=1}^Nk_m(z)
\]
is analytic in $\C_+$ and $v(z)=\Im w(z)>0$, $z\in\C_+$. Moreover, $w$ is continuous in
$\ol{\C}_+$.  Using this fact and Theorem \ref{T41} we deduce that
$v$ has an harmonic extension from $\C_+$ into
$\O$ by $v(z)=v(\ol z),\ \ z\in \C_\pm$.
Moreover, $w$ has an analytic extension from $\C_+$ into
$\C_+\cup \C_-\cup g_n$ since $u=$ const on the gap $g_n$, and we have identities \er{412}.

 We show that $v$ is subharmonic in $\C$.
Since $v$ is harmonic  in $\C_+\cup \C_-\cup g$ we have to show
that $v$ is subharmonic at $z_0\in \R\sm g$. Note that $v$ is
continuous in $\C$.

Introduce the sets
$$
\a_s(z)=\{m: \D_m(z)\in [-1,1]\},\ \ \a_g(z)=\{m: \D_m(z)\in
\R\sm [-1,1]\},\ \ \ \
$$
and $ \a_c(z)=\{m: \D_m(z)\notin \R\}$. We have the decomposition
$$
v=v_s+v_c+v_g, \ \  \ \ v_j=\sum _{m\in
\a_j(z_0)}{q_m\/N},\  \ j=s,c,g.
$$
 We have 2 cases: the first case $z_0\notin \cB_\D$. 

Consider $v_g$.  We take the interval
$Y=(z_0-\ve,z_0+\ve)$ for some $\ve>0$ such that $Y\cap
\cB_\D=\es$. Due to Theorem \ref{T41} each $q_m(\cdot), m\in
\a_g(z_0)$ has an harmonic extension from $\mR_0^+$ into
$\mR_0^+\cup\mR_0^-\cup Y$ by $q_m(z)=q_m(\ol z),\ \ z\in
\mR_0^\pm$ for sufficiently small $\ve$. Thus $v_g$ has an
harmonic extension from $\mR_0^+$ into
$\mR_0^+\cup \mR_0^-\cup Y$ by $v_g(z)=v_g(\ol z),\ \ z\in
\mR_0^\pm$.

Consider $v_c$. Due to Theorem \ref{T41} we deduce that
$$
v_c=\sum _{m\in \a_c(z_0)}q_m(z)={1\/2}\sum _{m\in
\a_c(z_0)}q_m(z)+q_{\wt m}(z),\
$$
where $k_{\wt m}(z), \wt m\in \a_c, $ is the quasimomentum such
that $\D_{\wt m}(z)=\ol\D_m(z), z\in Y$. Due to Theorem \ref{T41}
each $v_m=q_m(\cdot)+q_{\wt m}, m\in \a_g(z_0)$ has an harmonic
extension from $\mR_0^+$ into $\mR_0^+\cup\mR_0^-\cup Y$ by
$v_m(z)=v_m(\ol z),\ \ z\in \mR_0^\pm$ for sufficiently small
$\ve$. Thus $v_c$ has has an harmonic
extension from $\mR_0^+$ into $\mR_0^+\cup\mR_0^-\cup Y$ by
$v_c(z)=v_c(\ol z),\ \ z\in \mR_0^\pm$ .

Hence due to $v_s(z_0)=0$ we obtain
$
v(z_0)\le {1\/2\pi }\int_0^{2\pi}v(z_0+re^{i\f})d\f,
$
for any small $r>0$, 
and thus $v$ is subharmonic in in some neighborhood of  $z_0$.

The second case $z_0\in\cB_\D$. Define two intervals
$Y_-=(z_0-2\ve,z_0), Y_+=(z_0, z_0+2\ve),$ for some $\ve>0$ such
that $Y_\pm\cap \cB_\D=\es$ and each $q_m$ is harmonic in $\{|z-z_0|<\ve\}\cap \C_+$. Define the functions
$$
v=v_1+v_0,\ \ 
v_0(z)=\sum _{m\in \a_0}q_m(z),\ \ \ \a_0=\a_s(z_0+\ve)\cup \a_s(z_0+\ve).
$$
Note that $v_0(z_0)=0$. 
The above arguments show that the function 
$v_1$
has a harmonic extension from $\mR_0^+$ into $\mR_0^+\cup\mR_0^-\cup Y_-\cup
Y_+$ by $v_1(z)=v_1(\ol z)$ for sufficiently small $\ve>0$.
Then the function $v_1$ has an harmonic extension from $\mR_0^+$
into $\mR_0^+\cup\mR_0^-\cup (z_0-2\ve, z_0+2\ve)$ by
$v_1(z)=v_1(\ol z)$ for sufficiently small $\ve$, since the
function $v_1$ is bounded in the disk $B(z_0,t)=\{z: |z-z_0|<\ve\}$
for some small $\ve>0$. Thus we obtain
$
v(z_0)\le {1\/2\pi }\int_0^{2\pi}v(z_0+re^{i\f})d\f
$
for any small $r>0$, since $v_0(z_0)=0$.
Hence the function $v$ is subharmonic in $\C$.

Asymptotics \er{41}, \er{42} give \er{412},\er{413}.
Moreover, \er{413} yields $\int_{\R}{v(t)\/1+t^2}dt<\iy$ and
$\lim\sup_{z\to\iy}{v(z)\/|z|}=1$, i.,e., $v\in \cS\cC$.

The function $v(z)-y$ is harmonic in $\C_+$, it is positive on the
real line and $v(z)-y=o(1)$ as $|z|\to \iy, z\in \ol\C_+$. Then
$v(z)-y$ is positive on $\C_+$. Using \er{412} and the Nevanlinna
Theorem we deduce that $\int_{\R}v(t)dt<\iy$. Then $v\in
\cS\cK_0^+$. $\BBox$

\section {The conformal mappings}
\setcounter{equation}{0}

 We need the result from [KK1].

\begin{theorem} \lb{T51} Let
$v\in \cS\cK_0^+$ and $v\neq const$. Then $w:\C_+\to w(\C_+)=W(h)$
is a conformal mapping for some $h\in C_{us}, h\ge 0$. Moreover, the
following asymptotics, estimates and identities are fulfilled:
\[
\lb{51} 
w(z)=z-{Q_0+o(1)\/z},\ \ \ \ |z|\to \iy,\ \ y\ge r|x|,
\ \ {\rm for \ any}\ r>0,
\]
\[
\lb{52}  I_0^D+\mP_0=Q_0, \ \ \ \ \ \
\sup_{x\in\R}v^2(x)\le 2Q_0.
\]
If in addition, $Q_2<\iy$, then $Q_2=I_1^D+\cP_2$, where
$\cP_n, I_n^D, Q_n, n\ge 0$ are given by \er{dQSI}.
\end{theorem}

\no {\bf Proof of Theorem \ref{T3}.} By Theorem \ref{T42}, the
function $v={1\/N}\sum_1^{N}q_m\in \cS\cK_0^+$. Then Theorem \ref{T51}
gives that $w:\C_+\to w(\C_+)=W(h)$ is a conformal mapping for
some $h\in C_{us}$.  

The function $M(iy),y>0$ is real, then $M(z)=\ol M(-\ol z), z\in \ol\C_+$. This yields that the set $\{\t_m(z)\}_1^N=\{\t_m(-\ol z)\}_1^N,z\in \ol\C_+$, which gives $v(z)=v(-\ol z), z\in \ol\C_+$.
Thus $v(x)=v(-x)$ for all $x\in \R$ and the identities
$$
w(-\ol z)=-\ol z+{1\/\pi }\int_\R{v(t)dt\/t+\ol z}=-\ol z-{1\/\pi }\int_\R{v(s)dt\/s-\ol z}=-\ol w(z), z\in \C_+
$$
 give $-w_0(-\ol z)=\ol w_0(z), z\in \C_+$.

Using \er{1asD} and $q_m\ge 0 $ we obtain  $q_m(iy)=y+o(1)$ and
\[
\det L(z)=\prod_1^N
\D_m^2(z)={(1+e^{-4Ny}O(1))\/4^N}e^{-2i\sum_1^N k_m(z)}
={(1+e^{-4Ny}O(1))\/4^N}e^{-2iNw(z)}
\]
as $z=iy, y\to \iy$  and \er{318} yields
\[
\det L(z) =(\cos^{2N}
z)\exp \rt(i{\Tr V^0\/z}+{i\|V\|^2+o(1)\/4z^3}\rt).
\]
Thus due to $\Re w(iy)=0$ we get
\[
w(z)=z-{\Tr V^0\/2Nz}-{\|V\|^2+o(1)\/8Nz^3}, \ \ \
z=iy,\ y\to \iy.
\]
 Then by the Nevanlinna Theorem and Theorem \ref{T51}, we have
asymptotics \er{T3-2} and identities \er{T3-3}, \er{T3-4}.

Estimate \er{52} gives $v\big|_\R\le\sqrt{2Q_0}$.
If $z\in \s_{(N)}$, then $v(z)=0$, since $\D_m(z)\in[-1,1]$
 and $q_m(z)=0$ for all $m=1,..,N$. 
If $z\in \s_{(1)}\cup g$, then  $\D_m(z)\notin [-1,1]$ for some $m=1,..,N$. Thus $q_m(z)>0$ and $v(z)>0$,  
which gives \er{T3-5}. \BBox

\no {\bf Proof of Theorem \ref{T5}.} i) We show \er{T5-1}. Using
$w(z)=z+{1\/\pi}\int_{\R}{v(t)dt\/(t-x)}, z\in \C_+$ and  \er{T3-5}  we obtain
$w'(z)=1+{1\/\pi}\int_{\R}{v(t)dt\/(t-x)^2}>1,\ \ z\in \s_{(N)}$,
which yields $u_x'(z)>1,z\in \s_{(N)}$.
We used the fact that due to \er{T2-3} $u_x'(z)=1$ for some $z\in \s_{(N)}$ iff $V=0$.

Let $x\in \s_{(1)}$.  Then some branch $\D_m(x)\in (-1,1)$ for all $x\in Y\ss\s_{(1)}$ for some small interval  $Y=(\a,\b)$.  We have $\D_m(x)=\cos
k_m(x)$. Thus we get $k_m'(x)=-{\D_m'(x)\/\sin
k_m(x)}\ne 0$. 
Hence we get $u_m'(x)=k_m'(x)>0$, since $q_m(x)=0$ and \er{49}
yields $u_x(x)\ge 0$ for $x\in Y$, which gives \er{T5-1}.

We show \er{T5-2}. 
 Note that $v(z)>0,z\in g$, otherwise we have not a gap. Using \er{49} we obtain $ -{\pa^2 v(z)\/\pa x^2}
={\pa^2 v(z)\/\pa y^2}=\int_{\R}{d\m_v(t)\/(t-x)^2}>0,\ \ \ z=x\in
g_n, $ which yields $v_{xx}''(z)<0, z\in g_n$. Consider the function $u(z), z\in g_n$.
Theorem \ref{T41} yields
$
N\Re w(z)=\sum_1^N \Re k_m=\sum_1^N \pi n_m=\pi N_{n}, 
$
which gives \er{T5-2}

We recall the result from [KK2]. Let a function $f$ be harmonic
and positive in the domain $\C\sm [-a,a],\ a>0$ and
$f(iy)=y(1+o(1))$ as $y\to\iy$. Assume $f(z)=f(\ol z),\ z\in
\C\sm [-a, a]$ and $f\in C(\ol{\C_+})$. Then
$$
f(x)=\sqrt {a^2-x^2}\lt(1+{1\/\pi}\int_{\R\sm I} {f(t)dt\/
|t-x||t^2-a^2|^{1\/2}}\rt),\ \  x\in [-a,a].
$$
Hence the last identity and properties of $v$ yield \er{T5-3}
and the estimate $v(z)\ge v_n^0(z)=|(z-z_n^-)(z_n^+-z)|^{1\/2},\ \ \
z\in g_n=(z_n^-,z_n^+)$. Consider the case $g_n\ss \R_+$ (the proof for the cases $g_n\ss \R_-$ is similar).
Defining $z_n^0={z_n^++z_n^-\/2}, r={|g_n|\/2}$
and using $(z_0+x)^{2p}+(z_0-x)^{2p}\ge 2z_0^{2p}$ for all $p\ge 0$ we have 
$$
{1\/\pi}\int_{g_n}t^{2p}v(t)dt\ge
\int_{g_n}{t^{2p}v_n^0(t)dt\/\pi}\ge {1\/\pi}\int_0^r\sqrt
{r^2-x^2}((z_n^0+x)^{2p}+(z_n^0-x)^{2p})dx\ge
{r^2\/2}(z_n^0)^{2p}.
$$
 If $p=0$, then we get $Q_0={1\/\pi}\int_\R v(t)dt\ge {1\/8}
\sum |g_n|^2$, which yields the first estimate in \er{T5-4}.

If $p=1$, then using $\g_n=\l_n^+-\l_n^-\ge 0,\ \
\l_n^\pm=(z_n^\pm)^2$  we get
$$
Q_2={1\/\pi}\int_\R t^2v(t)dt\ge {1\/8} \sum_{n\in \Z}
|g_n|^2{(z_n^++z_n^-)^2\/4}={1\/16} \sum_{n\ge 0} |\g_n|^2,
$$
which yields the  last inequality in \er{T5-4}. 

Recall estimates from [K10].
 The conformal mapping $w:\C_+\to W(h)$ for the case $w(z)=z+o(1)$
as $|z|\to \iy$ and $h(u)=0, u\ne {\pi\/N} \Z$ and
$\{h({\pi\/N}n)\}_{-\iy}^\iy\in \ell_1^2$ was studied in [K7].
For some absolut constant $C_0$ the following estimate was obtained:
$Q_2\le C_0G^2(1+G^{1\/3})$, which together with
$Q_2={\|V\|^2\/8N}$ gives the estimate  \er{T5-5}.
\BBox

\section {Appendix}
\setcounter{equation}{0}

\begin{lemma} \lb{TA1}
Let function $F(s,t),R(t,s,u), t,s,u\in (0,1)$
satisfy : 

i) $F(\cdot,\cdot)\in L^2((0,1)^2)$ and $F(t,t)\in L^1(0,1)$,

ii) $R(\cdot,\cdot,\cdot)\in L^2((0,1)^3)$ and $ R(t,t,s),R(t,s,t),R(s,t,t)\in L^1((0,1)^2)$ .
Then
\[
\lb{61} f(z)=\int_0^1dt\int_0^t\cos z(1-2t+2s)F(t,s)ds ={i\cos
z\/2z}\lt(\int_0^1F(t,t)dt+o(1)\rt)
\]
\[
\lb{62} f^+(z)=\int_0^1dt\int_0^t
e^{-i2z(t-s)}F(t,s)ds={o(e^{2y})\/|z|},
\]
\[
\lb{63} \int_0^1dt\int_0^tds\int_0^s \sin
z(1-2\z)R(t,s,u)du=o(e^{y}),
\]
as $r|x|<y\to \iy$ for any fixed $r>0$, where $\z$ is one of functions: $s-u, t-u $ or $t-s$.
\end{lemma}

\no {\it Proof.} We have $f(z)=e^{iz}f^+(z)+e^{-iz}f^-(z)$,
where
\[
\lb{6a} f^-(z)\ev{1\/2}\int_0^1dt\int_0^te^{i2z(t-s)}F(t,s)ds
={1\/4}\int_0^1ds\int_0^1e^{i2z|t-s|}F(t,s)ds
\]
where $F(t,s)=F(s,t), t,s\in (0,1)$.  Substituting the identities
$$
{1\/ 2\pi }\int_\R  e^{ipk}{dk\/k^2-4z^2}={i\/4z}e^{i2z|p|},\ \
p\in \R,\ z\in \C_+,\ \ \ \ 
\hat F(r,k)\ev{1\/ 2\pi }\iint_{[0,1]^2}
e^{-irt-iks}F(t,s)dtds, 
$$
into \er{6a} we obtain
$$
f^-={z\/i}\int_\R {\hat F(-k,k)dk\/k^2-4z^2}
={1\/i4z}\int_\R \rt(-1+{k^2\/k^2-4z^2}\rt)\hat F(-k,k)dk=
{-1\/i4z}\rt(\int_0^1 F(t,t)dt+o(1)\rt).
$$
Consider $f^+$. Let $\|F\|^2=\int_0^1dt\int_0^t |F(t,s)|^2ds$  We have
\[
\lb{6b} |f^+(z)|^2\le g(z)\|F\|^2,\ \  \ \ g(z)=\int_0^1dt\int_0^t
e^{4y(t-s)}ds\le \int_0^1 {e^{4yt}\/4y}dt\le {e^{4y}\/(4y)^2}.
\]
Let $F_0$ be a smooth function such that $\|F-F_0\|\le \ve$ for
some small $\ve>0$. Define the function
$f_0^+(z)=\int_0^1dt\int_0^t e^{-i2z(t-s)}F_0(t,s)ds$. Using
\er{6b} we obtain
\[
|f^+(z)|\le |f_0^+(z)|+|f^+(z)-f_0^+(z)|\le |f_0^+(z)|+
\|F-F_0\|{e^{2y}\/(4y)}
\]
and the integration by parts yields $f_0^+(z)=O({e^{2y}\/y^2})$.
Thus we obtain \er{62}.

The proof for \er{63} is similar. $\BBox$

\begin{lemma} \lb{TA3}
The function $f(z)=\log |\x(z)|, z\in \C\sm [-1,1]$ is subharmonic
and continuous in $\C$. Moreover, for some absolute constant $C$
the following estimate is fulfilled:
\[
\lb{68} |f(z)-f(z_0)|\le C\ve^{1\/2},\ \ \ if \ \ |z-z_0|\le \ve
\max\{2,|z_0|\}, \ \ 0\le\ve\le {1\/8}, \ \ z,z_0\in \C.
\]
\end{lemma}
\no {\it Proof.} The function $f(z)=\log |\x(z)|0>, z\in \C\sm
[-1,1]$ is harmonic and $f(z)=0, z\in [-1,1]$. Then $f$ is
subharmonic in $\C$.

1. Let $|z_0|\le 2$. Consider the case $|z_0-1|\le 1/2$.
The proof of other cases is similar.
Firstly, let  $z,z_0\in \ol \C_+$. Then
$$
f(z)-f(z_0)=\Re \int_{z_0}^z{dt\/\sqrt{t^2-1}}=\Re
\int_{z_0-1}^{z-1}{ds\/\sqrt{s(2+s)}},\ \ \ \ \
{1\/\sqrt{s(2+s)}}={1\/\sqrt{s}}+H(s),
$$
where $H(\cdot)$ is analytic and bounded in the disk $B(0,1)$.
Then
$$
f(z)-f(z_0)=\Re 2\sqrt s \rt|_{z_0-1}^{z-1}+O(\ve )
$$
which yields \er{68}.

Consider the case $z\in\C_-,z_0\in \C_+$.
The identity $f(z)=f(\ol z)$ gives $f(z)-f(z_0)=f(\ol z)-f(z_0),
\ol z, z_o\in \C_-$. Then the simple estimate
$|\ol z-z_0|\le |z-z_0|$  and the case $z,z_0\in \ol \C_+$
imply \er{68}.

Let $|z_0|> 2$. Let $a_0=z_0/|z_0|,a=z/|z_0|$. Note that
$|a-a_0|\le \ve$.
$$
f(z)-f(z_0)=\Re \int_{z_0}^z{dt\/\sqrt{t^2-1}}=\Re
\int_{a_0}^aF(s)ds,\ \ \ F(s)=\rt(s^2-{1\/|z_0|^2}\rt)^{-{1\/2}},
\ \ s=t/|z_0|.
$$
 We obtain
$$
\rt|s^2-{1\/|z_0|^2}\rt|\ge |s|^2-{1\/4}\ge (1-\ve)-{1\/4}\ge
{1\/2},\ \ \ {\rm for} \ |s-a_0|\le \ve.
$$
Thus the function $F(s)=\rt(s^2-{1\/|z_0|^2}\rt)^{-{1\/2}}$ is
analytic and bounded in $s, |a_0-s|\le \ve$. Thus we deduce that $
|f(z)-f(z_0)|\le \ve C $ for some absolute constant $C$. $\BBox$

 \no {\bf Acknowledgments.} The authors would like
to thank Horst Hohberger for the Figure 2. 
Dmitry Chelkak was partly supported by grants VNP Minobrazovaniya 3.1--4733, RFFR 03--01--00377 and NSh--2266.2003.1. 
Evgeny Korotyaev was partly supported by DFG project BR691/23-1.

\no {\bf References}

\no [AS] Abramowitz, M. and Stegun, A., eds.:
Handbook of Mathematical Functions. N.Y.: Dover
Publications Inc., , 1992

\no [Ah] Akhiezer, N.: The classical moment problem and some
related questions in analysis. Hafner Publishing Co., New York
1965

\no [BBK] Badanin, A; Br\"uning, J; Korotyaev, E. The Lyapunov
function for Schr\"odinger operator  with periodic  $2\ts 2$
matrix potential, preprint 2005

\no [BK] Badanin, A; Korotyaev, E.
Spectral asymptotics  for periodic forth order operators, 
will be published in  Int. Math. Res. Not.

\no [C1] R. Carlson. An inverse problem for the matrix
Schr\"odinger equation. J. Math. Anal. Appl.
267 (2002), no. 2, 564--575.

\no [C2] R. Carlson. Eigenvalue estimates and trace formulas
for the matrix Hill's equation. J. Differential Equations 167
(2000), no. 1, 211--244.

\no [C3] R. Carlson. Compactness of Floquet isospectral sets for
the matrix Hill's equation. Proc. Amer. Math. Soc. 128 (2000),
 no. 10, 2933--2941.

\no [C4] Carlson, R. A spectral transform for the matrix Hill's equation. Rocky Mountain J. Math. 34 (2004), no. 3, 869--895.

\no [CG] Clark, S., Gesztesy, F.: Weyl-Titchmarsh $M$-function asymptotics, local uniqueness results, trace formulas, and Borg-type theorems for Dirac operators. Trans. Amer. Math. Soc. 354 (2002), no. 9, 3475--3534

\no [CG1] S. Clark; F. Gesztesy. Weyl-Titchmarsh $M$-function
asymptotics for matrix-valued Schr\"odinger operators.
Proc. London Math. Soc. (3) 82 (2001), no. 3, 701--724.

\no [CHGL] Clark S., Holden H., Gesztesy, F., Levitan, B.:
Borg-type theorem for matrix-valued Schr\"odinger and Dirac
operators, J. Diff. Eqs. 167(2000), 181-210

\no [CS]  Craig, W.; Simon, B. Subharmonicity of the Lyaponov index.
Duke Math. J.  50  (1983),  no. 2, 551--560.

\no [DS] Dunford, N. and Schwartz, J.: Linear Operators Part II:
Spectral Theory, Interscience, New York, 1988.

\no [Fo]  Forster O.: Lectures on Riemann surfaces. 
Graduate Texts in Mathematics, 81. Springer -Verlag, New York, 1991

\no  [GT1] J. Garnett, E. Trubowitz: Gaps and bands of one dimensional  periodic Schr\"odinger operators. Comment. Math. Helv. 59, 258-312 (1984)

\no  [GT2] J. Garnett, E. Trubowitz: Gaps and bands of one dimensional
 periodic Schr\"odinger operators II. Comment. Math. Helv. 62, 18-37 (1987).

\no [GL] Gel'fand I., Lidskii, V.: On the structure of the regions
of stability of linear canonical
   systems of differential equations with periodic coefficients.
   (Russian)
   Uspehi Mat. Nauk (N.S.) 10 (1955), no. 1(63), 3--40.

\no [Ge] Gel'fand, I.:
 Expansion in characteristic functions of an equation with
 periodic coefficients. (Russian) Doklady Akad. Nauk SSSR (N.S.) 73, (1950). 1117--1120.

\no [GKM] Gesztesy, F., Kiselev A.; Makarov, K. Uniqueness
results for matrix-valued Schrodinger, Jacobi, and Dirac-type
operators. Math. Nachr. 239/240 (2002), 103--145.

\no  [Ka] T. Kappeler:
Fibration of the phase space for the Korteveg-de-Vries equation.
Ann. Inst. Fourier (Grenoble), 41, 1, 539-575 (1991).

\no [KK] Kargaev P., Korotyaev E.: Inverse Problem for the Hill
Operator, the Direct Approach.  Invent. Math., 129(1997), no. 3,
567-593

\no [KK1] Kargaev P., Korotyaev E.:
Identities for the Dirichlet integral of subharmonic functions
from the Cartright class, 
 Complex Var. Theory Appl. 50 (2005), no. 1, 35--50.

\no [KK2] Kargaev P., Korotyaev E.:
Conformal mappings and subharmonic functions of Polya class

\no [Ka] Kato, T.: Perturbation theory for linear operators.
Springer-Verlag, Berlin, 1995.

\no [K1]  E. Korotyaev. The inverse problem for the Hill operator. I
Internat. Math. Res. Notices, 3(1997), 113--125

\no [K2]  E. Korotyaev. Inverse problem and the trace formula for the Hill operator. II Math. Z. 231(1999), no. 2,    345--368

\no [K3]  E. Korotyaev.  Inverse problem for periodic "weighted" operators,  J. Funct. Anal. 170(2000), no. 1, 188--218

\no [K4] E. Korotyaev.  Marchenko-Ostrovki mapping for periodic
Zakharov-Shabat systems,   J. Differential Equations, 175(2001), no. 2, 244--274

\no [K5] E. Korotyaev.  Inverse Problem and Estimates for Periodic
Zakharov-Shabat systems, J. Reine Angew. Math. 583(2005), 87-115

\no [K6] E. Korotyaev. Characterization of the spectrum of Schr\"odinger operators with periodic distributions. Int. Math. Res. Not.  (2003) no. 37, 2019--2031

\no [K7] Korotyaev, E.: The estimates of periodic potentials in
terms of effective masses. Comm. Math. Phys.
   183 (1997), no. 2, 383--400.
   
   \no [K8] Korotyaev, E.: Metric properties of conformal mappings on the complex plane with parallel slits, Internat. Math. Res.
Notices, 10(1996), 493--503

\no [K9]  Korotyaev, E.:  Estimates for the Hill operator. I J. Differential Equations 162 (2000), no. 1, 1--26

\no [K10]  Korotyaev, E.:   Estimates of periodic potentials in terms of gap lengths, Comm. Math. Phys. 197 (1998), no. 3, 521--526

\no [K11]  Korotyaev, E.:   Estimates for the Hill operator, II, will be published in J. Differential Equations,

\no [Kr] Krein, M.: The basic propositions of the theory of
$\lambda$-zones of stability of a canonical system of linear
differential equations with periodic coefficients. In memory of A.
A. Andronov, pp. 413--498. Izdat. Akad. Nauk SSSR, Moscow, 1955.

\no [Ly] Lyapunov, A.: The general problem of stability of motion,
2 nd ed. Gl. Red. Obschetekh. Lit., Leningrad, Moscow, 1935;
reprint Ann. Math. Studies, no. 17, Prinston Univ. Press,
Prinston, N.J. 1947

\no [MV] Maksudov, F.; Veliev, O.  Spectral analysis of differential operators with periodic matrix coefficients. (Russian) Differentsial'nye Uravneniya 25 (1989), no. 3, 400--409, 547; translation in Differential Equations 25 (1989), no. 3, 271--277

\no [MO1] Marchenko V., Ostrovski I.: A characterization of the
spectrum of the Hill operator. Math. USSR Sb. 26, 493-554 (1975).

\no [MO2] Marchenko V., Ostrovski I.: Approximation of periodic by
finite-zone potentials. Selecta Math. Sovietica. 1987, 6, No 2, 101-136.

\no [Mi1] Misura T. Properties of the spectra of periodic and anti-periodic
boundary value problems generated by Dirac operators. I,II, Theor. Funktsii Funktsional. Anal. i Prilozhen, (Russian), 30 (1978), 90-101; 31 (1979), 102-109

\no [Mi2] Misura T. Finite-zone Dirac operators.  Theor. Funktsii
Funktsional. Anal. i Prilozhen, (Russian), 33 (1980), 107-11.

\no [PT] P\"oshel,J., Trubowitz, E.: Inverse spectral theory. Pure
and Applied Mathematics, 130. Academic Press, Inc., Boston, MA,
1987. 192 pp.

\no [RS] M. Reed ; B. Simon. Methods of modern mathematical physics. IV. Analysis of operators. Academic Press, New York-London, 1978

\no [Sh] Shen, Chao-Liang: Some eigenvalue problems for the vectorial Hill's equation. Inverse Problems 16 (2000), no. 3, 749--783.

\no [Sp] Springer, G.: Introduction to Riemann surfaces. Addison-Wesley Publishing Company, Inc., Reading, Mass. 1957 

\no [Ti] Titchmarsh, E.: Eigenfunction expansions
 associated with second-order differential equations 2, Oxford: Clarendon Press 1958

\no [YS] Yakubovich, V., Starzhinskii, V.: Linear differential
equations with periodic coefficients. 1, 2.
   Halsted Press [John Wiley \& Sons] New York-Toronto,
   1975. Vol. 1, Vol. 2

\end{document}